\documentclass{svmult}
\usepackage{amsmath}
\usepackage{amssymb}

\makeindex             

\begin{document}

\title*{Canonical Relativized Cylindric Set Algebras and Weak
  Associativity} \titlerunning{Canonical Relativized Cylindric Set
  Algebras and Weak Associativity} \author{Roger D.  Maddux}
\institute{Roger D. Maddux \at Department of Mathematics, Iowa State
  University, Ames, Iowa 50011-2066, USA, \email{maddux@iastate.edu}}

\maketitle

\abstract{Canonical relativized cylindric set algebras are used to
  sharpen the relative representation theorem for weakly associative
  relation algebras, that every complete atomic weakly associative
  relation algebra is isomorphic with the relativization of a set
  relation algebra to a symmetric and reflexive binary relation, by
  insuring that the atoms of the set relation algebra and its
  relativization are orbits of single sequences under a group of
  permutations of the underlying set.  This sharpening of the relative
  representation theorem was first proved for the Resek-Thompson
  Theorem in 1989.}

\medskip\noindent{\bf Keywords\ } Cylindric set algebras $\cdot$ Set
relation algebras $\cdot$ Relativized cylindric algebras $\cdot$
Relativized relation algebras $\cdot$ Symmetric relation $\cdot$
Reflexive relation $\cdot$ Relativization
\date{July 15, 2018, revised October 16, 2019, corrected and revised
  February 16, 2021}

\let\SS=\S
\allowdisplaybreaks
\newtheorem{cor}{Corollary}

\def\con#1{\setbox13\hbox{$#1$}\ifdim\wd13<.9em\breve{#1}\else{\(#1\)}\breve{\ }\fi}
\def \halfthinspace{\relax\ifmmode\mskip.5\thinmuskip\relax\else\kern.8888em\fi}
\def \({\left(}
\def \){\right)}
\def \<{\langle}
\def \>{\rangle}
\def \minus{\setminus}
\def \lb{\left\{\,}
\def \rb{\,\right\}}
\def \gc{\mathfrak}
\def \ga{{\gamma}}
\def \Rs{\gc{Rs}\halfthinspace}
\def \Rc{\gc{Rc}\halfthinspace}
\def \Rl{\gc{Rl}\halfthinspace}
\def \Sb{\gc{Sb}\halfthinspace}
\def \Cm{\gc{Cm}\halfthinspace}
\def \Ra{\gc{Ra}\halfthinspace}
\def \Re{\gc{Re}}
\def \Tr{Tr} 
\def \Fd{Fd} 
\def \ato#1{{}^{\alpha}{#1}} 
\let \ism=\cong
\def \mi{^{-}}
\def \ie{{\it i.e.\/}}
\def \eg{{\it e.g.\/}}

\def\BA{\mathsf{BA}}
\def\NA{\mathsf{NA}}
\def\WA{\mathsf{WA}}
\def\SA{\mathsf{SA}}
\def\RA{\mathsf{RA}}
\def\At{\mathsf{At\,}}
\def\at{\mathsf{At\,}}
\def\id{1'}
\def\id{{1\kern-.08em\raise1.3ex\hbox{\rm,}\kern.08em}}
\def\di{{0\kern-.04em\raise1.3ex\hbox{,}\kern.04em}}
\def\rp{{;}}
\def\halfthinspace{\relax\ifmmode\mskip.5\thinmuskip\relax\else\kern.8888em\fi}
\let\hts=\halfthinspace
\def\rp{{\hts;\hts}}

\def\sub{\mathsf{s}}
\def\cyl{\mathsf{c}}
\def\diag{\mathsf{d}}
\def\CA{\mathsf{CA}}
\def\halfthinspace{\relax\ifmmode\mskip.5\thinmuskip\relax%
\else\kern.8888em\fi}
\def\Cs{\mathsf{Cs}}
\def\Ca{\gc{Ca}\halfthinspace}
\def\rep#1{R^{\gc\B}_{#1}}
\def\repp{\mathsf{F}}
\def\RL{\mathsf{Rl}\halfthinspace}
\def\etc{{\it etc.}}
\def\NCA{\mathsf{NCA}}
\def\Cyl{\mathsf{C}}
\def\Diag{\mathsf{D}}
\def\dom#1{\setbox13\hbox{$#1$}\ifdim\wd13<11pt{#1}^{\mathsf{d}}%
\else{\(#1\)}^{\mathsf{d}}\fi}
\def\rng#1{\setbox13\hbox{$#1$}\ifdim\wd13<11pt{#1}^{\mathsf{r}}%
\else{\(#1\)}^{\mathsf{r}}\fi}
\let\bp=\cdot
\def\dash{{\text{-}}}
\def\min#1{\overline{#1}}
\def\sigop#1{\setbox13\hbox{$#1$}\ifdim\wd13=0pt%
\z\rp{#1}\else\z\rp{#1}\fi}
\def\tauop#1{\setbox13\hbox{$#1$}\ifdim\wd13=0pt%
\con\z\rp{#1}\else\con\z\rp{#1}\fi}
\def\rhoop#1{\setbox13\hbox{$#1$}\ifdim\wd13=0pt%
\rho{#1}\else\rho\({#1}\)\fi}
\def\blank{\phantom{x}}
\def\atm{\at{\gc\A}}
\def\pow{\mathcal{P}}

\def\makecs#1#2{\makecsX {#1}#2,.}
\def\makecsX#1#2#3.{\onecs{#1}{#2}
\ifx#3,\let\next\eatit\else\let\next\makecsX\fi\next{#1}#3.}
\def\onecs#1#2{\expandafter\gdef\csname #2\endcsname%
{{\csname #1\endcsname {#2}}}}
\def\eatit#1#2.{\relax}
\makecs{}{abcdef hij  mnopqrstuvwxyzABCDEFGHIJKLMNOPQRSTUVWXYZ}
\def\g{\mathcal{A}}
\def\k{\kappa}
\def\l{\lambda}
\def\m{\mu}
\def\un{{}^3\U(\gc\B)}

\section{Introduction}\label{sect0}
My first meeting and mathematical encounter with Hajnal and Istv\'an
occurred in 1984 on a beach in Charleston, South Carolina, at a
meeting organized by Steve Comer. As the three of us came out of the
water, Hajnal asked me, ``How do you prove\dots'', and I answered,
``Read the proof in my dissertation, \emph{not} the published proof.''
Her question referred to a theorem about relative representation,
proved in my dissertation \cite[Theorem 5(10)]{M78} for the class
$\SA$ of semi-associative relation algebras, and published in
\cite[Theorem 5.20]{MR662049} for the strictly larger class $\WA$ of
weakly associative relation algebras.  The original proof in
\cite[pp.\ 76--83]{M78} is more direct, and can be illustrated, as I
did on the beach, by drawing triangles in the sand (\'a la Archimedes)
that share edges. On any given edge, labelled with an atom included in
the product of two other atoms, make a new triangle by drawing two new
edges labelled with the two atoms in the product, and continue {\it ad
  infinitum}. See the proof of Theorem \ref{(10)} in \SS\ref{sect4}.

Axioms \eqref{ra1}--\eqref{ra10} for $\WA$ are given in
\SS\ref{sect1}.  Elementary consequences \eqref{1'sym}--\eqref{1(39)}
of these axioms, used in \SS\SS\ref{sect1}--\ref{sect10}, are
relegated (along with complete proofs) to \SS\ref{sect11}.  Weak
associativity \eqref{ra4} is needed only for
\eqref{zx1880}--\eqref{1(39)}.  The closure of $\WA$ under the
formation of canonical extensions is discussed in \SS\ref{sect2}.  In
\SS\ref{sect3} the algebra $\Re(\U)$ of all binary relations on a set
$\U$ and the relativizations of its subalgebras are defined.
Corollary \ref{th1} shows relativizing $\Re(\U)$ to a symmetric and
reflexive relation yields a weakly associative relation algebra.
\SS\ref{sect4} presents a converse to this observation in Theorem
\ref{(10)}, whose proof (not previously published and quoted here from
\cite[pp.\ 76--83]{M78}) I recommended to Hajnal and Istv\'an: every
$\gc\A\in\WA$ is isomorphic to a subalgebra of
$\Rl_\E\big(\Re(\U)\big)$ for some symmetric and reflexive $\E$.
Theorem \ref{(10)} is used to characterize $\WA$ in Corollary
\ref{charWA} as the class of algebras isomorphic to subalgebras of
relativizations (to symmetric and reflexive relations) of the algebras
of binary relations on sets.  In \SS\SS\ref{sect5}--\ref{sect10} we
prove (in Theorem \ref{th7}) a sharpening of Corollary \ref{charWA}
for complete atomic algebras in $\WA$, namely, that every such algebra
is a relativization of a complete atomic set relation algebra whose
atoms are orbits of single sequences under a group of permutations of
the underlying set.  The proof uses results from \cite{MR987611},
where a similar sharpening of the Resek-Thompson Theorem \cite{AT} is
presented.

In \SS\ref{sect5} we recall from \cite{MR987611} the concepts of
suitable structure (of any dimension) and, for an arbitrary suitable
structure $\gc\B$, its canonical relativized cylindric set algebra
$\Rc\gc\B$.  In \SS\ref{sect6}, we construct from any atomic
$\gc\A\in\WA$ a 3-dimensional suitable structure $\gc\B(\gc\A)$.  In
\SS\ref{sect7} we define the complex algebra $\Cm\gc\B$ for any
suitable structure $\gc\B$. We show (in Lemmas \ref{NA3} and
\ref{MGR}) that for any atomic $\gc\A\in\WA$, $\Cm\gc\B(\gc\A)$ is a
non-commutative 3-dimensional cylindric algebra satisfying the
merry-go-round identities. The main result of \cite[Theorem
C]{MR987611} applies to such algebras, yielding an isomorphism between
$\Cm\gc\B(\gc\A)$ and $\Rc\gc\B(\gc\A)$ in Theorem \ref{th4}.  In
\SS\ref{sect8} we recall the definition of relation-algebraic reduct
$\Ra\gc\C$ of a non-commutative 3-dimensional cylindric algebra
$\gc\C$, and show, in Theorem \ref{th5}, that every complete atomic
$\gc\A\in\WA$ is isomorphic to $\Ra\Cm\gc\B(\gc\A)$.  The proof of
Theorem \ref{th6} in \SS\ref{sect9} combines the isomorphisms of
Theorem \ref{th4} and Theorem \ref{th5} with the definition of
canonical relativized cylindric set algebra from \SS\ref{sect5} to
show that every complete atomic $\gc\A\in\WA$ is isomorphic to the
relation-algebraic reduct of the relativization, to a ternary relation
$\V(\gc\B(\gc\A))$, of a complete atomic 3-dimensional cylindric set
algebra $\gc\C$, giving us the first of the following three
isomorphisms:
\begin{align*}
  \gc\A \cong\Ra\Rl_{\V(\gc\B(\gc\A))}\gc\C
  \cong\Rl_{\Cyl_2\V(\gc\B(\gc\A))}\Ra\gc\C \cong\Rl_\E\gc\A'.
\end{align*}
Most of the proof of Theorem \ref{th7} in \SS\ref{sect10} is devoted
to the second of the isomorphisms above. In this case, relativization
commutes with relation-algebraic reduct.  According to the link in
\cite{HMT85} between cylindric set algebras and relation set algebras,
the third isomorphism holds for some set relation algebra $\gc\A'$
isomorphic to $\Ra\gc\C$ and some symmetric and reflexive binary
relation $\E$, namely the image of the ternary relation
$\Cyl_2\V(\gc\B(\gc\A))$.  The elements of the base set of $\gc\A'$
are equivalence classes of certain finite sequences of triples of
atoms of $\gc\A$ alternating with members of $\{0,1,2\}$.
\section{Definition of $\WA$}\label{sect1}
$\WA$ is the class of {\bf weakly associative relation algebras},
algebras of the form $\gc\A=\<\A,+,\min\blank,\rp,\con\blank,\id\>$
satisfying equational axioms \eqref{ra1}--\eqref{ra10} below. These
equations are obtained from Tarski's classic axiomatization of the
class $\RA$ of relation algebras (see \cite[Definition
  2.1]{zbMATH06783012}) by replacing the {\bf associative law}
$(\x\rp\y)\rp\z=\x\rp(\y\rp\z)$ with the {\bf weak associative law}
\eqref{ra4}, accompanied by the definitions
\begin{equation*}
  \x\bp\y=\min{\min\x+\min\y},\qquad
  \di=\min\id,\qquad1=\id+\di,\qquad0=\min{\id+\di}.
\end{equation*}
Axioms \eqref{ra1}--\eqref{ra3} assert that $\<\A,+,\min\blank\,\>$ is
a Boolean algebra. If the weak associative law \eqref{ra4} were
replaced by the {\bf semi-associative law} $(\x\rp1)\rp1=\x\rp1,$ then
the axioms \eqref{ra1}--\eqref{ra10} would define the class $\SA$ of
semi-associative relation algebras.
\begin{align}
  \label{ra1}\x+\y&=\y+\x&&\text{$+$ is commutative}\\
  \label{ra2}\x+(\y+\z)&=(\x+\y)+\z&&\text{$+$ is associative}\\
  \label{ra3}\min{\min\x+\min\y}+\min{\min\x+\y}&=\x
  &&\text{Huntington's axiom}\\
  \label{ra4}((\x\bp\id)\rp1)\rp1&=(\x\bp\id)\rp1
  &&\text{weak associative law}\\
  \label{ra5}(\x+\y)\rp\z&=\x\rp\z+\y\rp\z&&\text{right additivity}\\
  \label{ra6}\x\rp\id&=\x&&\text{right identity law}\\
  \label{ra7}\con{\con\x}&=\x&&\text{involution}\\
  \label{ra8}\con{\x+\y}&=\con\x+\con\y
  &&\text{converse distributivity}\\
  \label{ra9}\con{\x\rp\y}&=\con\y\rp\con\x
  &&\text{converse-product rule}\\
  \label{ra10}\con\x\rp\min{\x\rp\y}+\min\y&=\min\y
  &&\text{Tarski/De~Morgan axiom}
\end{align}
All the consequences we need are proved in \SS\ref{sect11}. Among them
are equations that can form alternative axiom sets for $\WA$. The
axiomatization based on \cite[Theorem 1(22)]{M78} uses axioms
\eqref{ra1}--\eqref{ra4}, \eqref{ra6}, \eqref{ra7}, the definitions of
$\cdot$ and the constants, plus \eqref{x;0}, \eqref{0;x},
\eqref{leftid}, and \eqref{zx1766} from \SS\ref{sect11}.  Replacing
\eqref{ra4} with either the semi-associative law or the associative
law gives alternative axiomatizations for $\SA$ and $\RA$ as well.
\section{Canonical Extensions}\label{sect2}
A {\bf positive equation} is one that involves only $+$, $\bp$,
constants, $\rp$, and $\con\blank$, but not complementation.  Examples
include axioms \eqref{ra1}, \eqref{ra2}, \eqref{ra4}--\eqref{ra9}, the
alternative axioms listed in \SS\ref{sect1}, the semi-associative law,
and the associative law.  J\'onsson-Tarski \cite{MR0044502,MR0045086}
proved that every Boolean algebra with completely additive operators
has a canonical extension, and the extension satisfies the same
positive equations as the original algebra.  From this they conclude
that the canonical extension of a relation algebra is a relation
algebra.  For details see \cite{MR1935083}, \cite{MR2269199}, or
\cite{zbMATH06783011} (recommended).  The J\'onsson-Tarski results
apply to $\WA$ because $\rp$ and $\con\blank$ are completely additive
by \eqref{zx1699} and \eqref{ca} in \SS12.  Indeed, by \cite[Theorem
4.2]{MR662049}, every weakly associative relation algebra has a
canonical extension that is in $\WA$. (Since the associative laws are
positive, the same applies to $\SA$ and $\RA$.) The canonical
extension of an algebra is complete and atomic, but all we need in the
proof of Theorem \ref{(10)} is that every $\WA$ is a subalgebra of an
atomic $\WA$.
\section{Algebras of Binary Relations}\label{sect3}
For any set $\U$, the {\bf identity relation on} $\U$ is
\begin{align*}
  Id_\U=\{\<\u,\u\>:\u\in\U\}.
\end{align*}
The {\bf algebra of all binary relations on} $\U$ is
\begin{equation*}
  \Re(\U)=\<\pow(\U\times\U),\cup,{}_\U{\minus},\,|\,,{}^{-1},Id_\U\>
\end{equation*}
where $\pow$ is the power set operator, and for all
$\R,\S\subseteq\U\times\U$,
\begin{align*}
  \R\cup\S&=\{\<\u,\v\>:\<\u,\v\>\in\R\text{ or }\<\u,\v\>\in\S\},\\
  {}_\U{\minus}\R&=\{\<\u,\v\>:\u,\v\in\U\text{ and }\<\u,\v\>\notin\R\},\\
  \R|\S&=\{\<\u,\w\>:\exists\v(\<\u,\v\>\in\R
  \text{ and }\<\v,\w\>\in\S)\},\\
  \R^{-1}&=\{\<\u,\v\>:\<\v,\u\>\in\R\}.
\end{align*}
The {\bf field of} a binary relation $\E$ is
\begin{align*}
  \Fd(\E)=\{\u:\exists\v\big(\<\u,\v\>\in\E\text{ or }
  \<\v,\u\>\in\E\big)\}.
\end{align*}
A relation $\E$ is {\bf reflexive} if it contains the identity
relation on its field, \ie,
\begin{align*}
  Id_{\Fd(\E)}\cap(\E|\E^{-1}\cup\E^{-1}|\E)\subseteq\E,
\end{align*}
and {\bf symmetric} if $\E^{-1}=\E$.  If $\gc\A$ is a subalgebra of
$\Re(\U)$ and $\E\subseteq\U\times\U$ then the {\bf relativization of
  $\gc\A$ to} $\E$ is the algebra
\begin{align*}
  \Rl_\E\(\gc\A\)=\big\< \{\R\cap\E:\R\in\A\}, +, \min\blank, \rp,
  \con\blank, \id \big\>
\end{align*}
where, for all $\R,\S\in\A$,
\begin{align*}
  \R+\S&=(\R\cup\S)\cap\E,\\
  \min\R&=({}_\U{\minus}\R)\cap\E,\\
  \R\rp\S&=(\R|\S)\cap\E,\\
  \con\R&=\R^{-1}\cap\E,\\
  \id&=Id_\U\cap\E.
\end{align*}
The following direct consequence of the relevant definitions is also a
special case of \cite[Theorem 5.8(2)]{MR662049}, that the
relativization of a $\WA$ (such as $\Re(\U)$) to a symmetric-reflexive
element (\eg, any symmetric and reflexive relation) is in $\WA$.
\begin{cor}\label{th1}
  If $\U$ is a set, $\E\subseteq\U\times\U$, $\E$ is symmetric, and
  $\E$ is reflexive, then the relativization of $\Re(\U)$ to $\E$ is a
  complete atomic weakly associative relation algebra:
  $\Rl_\E\big(\Re(\U)\big)\in\WA$.
\end{cor}
\section{Characterizing $\WA$}\label{sect4}
In the proof of the following converse to Corollary \ref{th1}, a
relation $\E$ is constructed as a graph by adding a new vertex at each
step, along with two new edges connecting that vertex to the ends of a
previously selected edge. In the proof, Lemma \ref{A} asserts that
each step is possible, while Lemmas \ref{B} and \ref{C} say the
process can be completed to give the desired relation $\E$. The whole
construction is mediated by a labelling of the edges with atoms of an
atomic $\WA$, called a ``labelling system''. (The theorem was
originally stated and proved for $\SA$ but it applies to $\WA$ with no
changes.)
\begin{theorem}[{\rm\cite[Theorem 5(10)]{M78}}]\label{(10)}
  Let $\gc\A\in\WA$. Then there is a set $\U$ and a symmetric and
  reflexive binary relation $\E\subseteq\U\times\U$ such that $\gc\A$
  is isomorphic to a subalgebra $\gc\A'$ of the relativization to $\E$
  of the algebra of all binary relations on $\U$:
  \begin{equation*}
    \gc\A\cong\gc\A'\subseteq\Rl_\E\big(\Re(\U)\big).
  \end{equation*}
\end{theorem}

\proof \cite[pp.\ 76--83]{M78} We shall prove the theorem under the
additional assumption that $\gc\A$ is atomic. The original theorem
will then follow from the fact that every $\WA$ is a subalgebra of an
atomic $\WA$ (namely, its canonical extension).  Suppose a set $\U$
and a set $T\subseteq\U\times\atm\times\U$, where $\atm$ is the set of
atoms of $\gc\A$, have been constructed so that the following six
properties hold for all $a,b,c\in\atm$ and for all $u,v,w\in\U$:
\begin{enumerate}
\item[(i)] if $\<u,a,v\>,\<u,b,v\>\in T$, then $a=b$.
\item[(ii)] if $a\leq\id$ then there exists some $u\in\U$ such that
  $\<u,a,u\>\in T$,
\item[(iii)] $u=v$ iff there exists some $a\in\atm$ such that
  $a\leq\id$ and $\<u,a,v\>\in T$.
\item[(iv)] if $\<u,a,v\>\in T$ then $\<v,\con{a},u\>\in T$,
\item[(v)] if $\<u,a,v\>, \<v,b,w\>, \<u,c,w\>\in T$, then $a\rp b\geq
  c$,
\item[(vi)] if $\<u,a,v\>\in T$ and $a\leq b\rp c$, then there exists
  some $w\in\U$ such that $\<u,b,w\>,\<w,c,v\>\in T$.
\end{enumerate}
For all $x\in A$, set
\begin{equation*}
  \F(x)=\{\<u,v\>:u,v\in\U,\,\<u,a,v\>\in T\text{ for some $a\in\atm$ such
    that $a\leq x$}\},
\end{equation*}
and let $\E=\F(1)$.  We can then use properties (i)--(vi) to prove that
$\F$ is an isomorphism from $\gc\A$ into the algebra of all binary
relations contained in $\E$:
\begin{equation*}
  \F\in\text{Ism}\big(\gc\A,\Rl_\E\big(\Re(\U)\big)\big).
\end{equation*}
To prove the above, it is enough to show that conditions (vii)--(xii)
below hold for all $x,y\in A$.

\medskip
\par\noindent(vii) $\F(x+y)=\F(x)\cup \F(y)$.

\proof Suppose $\<u,v\>\in \F(x+y)$.  Hence $\<u,a,v\>\in T$ for some
$a\in\atm$ with $a\leq x+y$.  Then we get $a\leq x$ or $a\leq y$ since
$a$ is an atom, so either $\<u,v\>\in \F(x)$ or $\<u,v\>\in
\F(y)$. This shows $\F(x+y)\subseteq \F(x)\cup \F(y)$.  On the other
hand, if $\<u,v\>\in \F(x)$, then for some $a\in\atm$ we have
$\<u,a,v\>\in T$ and $a\leq x$, hence $a\leq x+y$, so $\<u,v\>\in
\F(x+y)$. Thus $\F(x)\subseteq \F(x+y)$, and similarly,
$\F(y)\subseteq \F(x+y)$.\endproof

\noindent(viii) $\F(\min\x)=\E\minus \F(x)$.

\proof Let $\<u,v\>\in \F(\min\x)$.  Then $\<u,a,v\>\in T$ and
$a\leq\min\x$ for some $a\in\atm$.  Suppose $\<u,v\>\in \F(x)$, \ie,
$\<u,b,v\>\in T$ for some $b\in\atm$.  Then $a=b$ by (i), and the
conditions $a\leq\min\x$ and $b\leq x$ are contradictory. Thus
$\<u,v\>\notin \F(x).$ From (vii) we get $\<u,v\>\in
\F(\min\x)\subseteq \F(\min\x+x)=\F(1)=\E$.  Thus $\<u,v\>\in \E\minus
\F(x)$.  Now let $\<u,v\>\in \E\minus \F(x)$. Then there is some
$a\in\atm$ such that $\<u,a,v\>\in T$ and for no $b\in\atm$ is it the
case that $\<u,b,v\>\in T$ and $b\leq x$.  In particular, since
$\<u,a,v\>\in T$, it is not the case that $a\leq x$.  This implies
that $a\leq\min\x$ since $a$ is an atom, and hence $\<u,v\>\in
\F(\min\x)$.\endproof

\noindent(ix) $\F$ is one-to-one.

\proof From (vii) and (viii) we conclude that $\F$ is a homomorphism
from the Boolean part of $\gc\A$ into the Boolean algebra of all
subrelations of $\E$:
\begin{equation*}
  \F\in\text{Hom}\big(Bl\,\gc\A,Bl\,\Rl_\E\big(\Re(\U)\big)\big).
\end{equation*}
We therefore need only show that $\F(x)\neq\emptyset$ whenever $0\neq
x\in A$. Suppose $0\neq x\in A$. Since $\gc\A$ is atomic, there is
some $a\in\atm$ with $a\leq x$. Then $a\rp\con{a}\cdot\id\in\atm$ by
\eqref{domrng} and \eqref{zx1890}, so by (ii) there is some $u\in\U$
such that $\<u,a\rp\con\a\cdot\id,u\>\in T$.  Now
$\a\rp\con\a\cdot\id\leq\a\rp\con\a$, so by (vi) there is some
$w\in\U$ such that $\<u,a,w\>,\<w,\con{a},u\>\in T$. But then
$\<u,w\>\in \F(x)$, so $\F(x)\neq\emptyset$.\endproof

\noindent(x) $\F(\id)=Id_\U$.

\proof Note that (x) is just a restatement of (iii).\endproof

\noindent(xi) $\F(\con{x})=\F(x)^{-1}$.

\proof Let $\<u,v\>\in \F(\con{x})$, so that $\<u,a,v\>\in T$ and
$a\leq\con{x}$ for some $a\in\atm$. Then $\<v,\con{a},u\>\in T$ by
(iv), $\con{a}\leq\con{\con{x}}=x$ by \eqref{ra7}, and
$\con{a}\in\atm$ by \eqref{1(36)}, so $\<v,u\>\in \F(x)$, and hence
$\<u,v\>\in \F(x)^{-1}$.  This shows that $\F(\con{x}) \subseteq
\F(x)^{-1}$ and hence also $\F(\con{x})^{-1}\subseteq\F(x)$.  Since
the latter formula holds for all $x\in A$, we can substitute $\con{x}$
for $x$ and obtain $\F(x)^{-1}=\F(\con{\con{x}})^{-1}\subseteq
\F(\con{x})$, which completes the proof of (xi).\endproof

\noindent(xii) $\F(x\rp y)=(\F(x)|\F(y))\cap \E$.

\proof Suppose $\<u,v\>\in \F(x\rp y)$, \ie, $\<u,a,v\>\in T$ and
$a\leq x\rp y$ for some $a\in\atm$.  Since $\gc\A$ is atomic and $\rp$
is completely additive by \eqref{ca},
\begin{equation*}
  a\leq x\rp y=\sum_{x\geq b\in\atm,\,\, y\geq c\in\atm}b\rp c,
\end{equation*} 
so there must be some $b,c\in\atm$ with $a\leq b\rp c$, $b\leq x$, and
$c\leq y$.  By (vi) there is some $w\in\U$ such that
$\<u,b,w\>,\<w,c,v\>\in T$.  From these facts we get $\<u,w\>\in
\F(x)$ and $\<w,v\>\in \F(y)$, and hence $\<u,v\>\in \F(x)|\F(y)$. We
also have $\<u,v\>\in \F(x\rp y)\subseteq \F(1)=\E$, so $\F(x\rp
y)\subseteq(\F(x)|\F(y))\cap \E$.

Now suppose $\<u,v\>\in(\F(x)|\F(y))\cap \E$. It follows that there
are $a,b,c\in\atm$ and $w\in\U$ such that $a\leq x$, $b\leq y$, and
$\<u,a,w\>,\<w,b,v\>,\<u,c,v\>\in T$. Then $a\rp b\leq x\rp y$ by
\eqref{rightmon} and \eqref{leftmon}, and $a\rp b\geq c$ by (v), so
$\<u,v\>\in \F(x\rp y)$. Thus $(\F(x)|\F(y))\cap \E\subseteq \F(x\rp
y)$, and the proof of (xii) is complete.\endproof

We have shown that
$\F\in\text{Ism}\big(\gc\A,\Rl_\E\big(\gc{Re}(\U)\big)\big)$. In
addition, using (xi), (vii), (x), and \eqref{1sym} we get
$\E^{-1}=\F(1)^{-1}=\F(\con1)=\F(1)=\E$ and $Id_\U=\F(\id)\subseteq
\F(1)=\E$, so $\E\subseteq\U\times\U$ is a symmetric and reflexive
relation. It follows that the theorem will be proved if we can find
sets $\U$ and $T$ satisfying (i)--(vi).  It is easy to get sets $\U$
and $T$ which satisfy (i)--(v). For example, let
$\U=\I^{(\gc\A)}=\{a:\id\geq\a\in\atm\}$ and
$\T=\{\<\a,\a,\a\>:\a\in\I^{(\gc\A)}\}.$ Then (i)--(iii) obviously
hold. To verify (iv) and (v) we only need to know that
$\a=\con\a=\a\rp\a$ whenever $a\in I^{(\gc\A)}$ (use \eqref{zx1850}
and \eqref{zx1849}).

We shall call $\<\U,T\>$ a {\bf labelling system} if
$T\subseteq\U\times\atm\times\U$ and properties (i)--(v) hold. (This
terminology was chosen because it is convenient to think of each
triple $\<u,a,v\>\in T$ as an ordered pair $\<u,v\>$ which has the
atom $a$ as a label.) Thus
\begin{equation*}
  \left<I^{(\gc\A)},\{\<a,a,a\>:a\in\I^{(\gc\A)}\}\right>
\end{equation*} 
is a labelling system. If property (vi) also holds, then we say the
labelling system $\<\U,T\>$ is {\bf complete}.  If $\<\U,T\>$ and
$\<\U',T'\>$ are labelling systems, then $\<\U',T'\>$ {\bf extends}
$\<\U,T\>$ if $\U\subseteq\U'$ and $T\subseteq T'$.  We wish to show
every labelling system can be extended to a complete labelling system.
A {\bf flaw} in a labelling system $\<\U,T\>$ is a quintuple
$\<u,a,v,b,c\>$ where $a,b,c\in\atm$, $u,v\in\U$, $\<u,a,v\>\in T$,
$a\leq b\rp c$, and for every $w\in\U$ either $\<u,b,w\>\notin T$ or
$\<w,c,v\>\notin T$.  Thus $\<\U,T\>$ is complete iff it has no flaws.
\begin{lemma}\label{A} If $\<\U,T\>$ is a labelling system and
  $\<u,a,v,b,c\>$ is a flaw in $\<\U,T\>$, then there is a labelling
  system $\<\U',T'\>$ which extends $\<\U,T\>$, such that
  $\<u,a,v,b,c\>$ is not a flaw in $\<\U',T'\>$.
\end{lemma}
\proof Choose any $\w$ such that $\w\notin\U$. Set $\U'=\U\cup\{\w\}$,
and
\begin{equation*}
  T'=T\cup\{\<\u,\b,\w\>,\<\w,\c,\v\>,\<\w,\con\b,\u\>,
  \<\v,\con\c,\w\>, \<\w,\rng\b,\w\>\},
\end{equation*}
where $\rng\b$ is defined in \eqref{domrng}, and
$\con\b,\con\c,\rng\b\in\atm$ by \eqref{1(36)} and \eqref{zx1890}.
Thus $T'\subseteq\U'\times\atm\times\U'$.  It is easy to check that
(i) and (iii) hold in $\<\U',\T'\>$.  (ii) holds for $\<\U',\T'\>$
simply because it holds for $\<\U,\T\>$.  (iv) is easily verified by
using \eqref{ra7} and \eqref{zx1849}.  It is obvious from the
construction of $\U'$ and $\T'$ that $\<\u,\a,\v,\b,\c\>$ is not a
flaw in $\<\U',\T'\>$. It remains to check (v).  Suppose
$\u_0,\u_1,\u_2\in\U'$, $\x_0,\x_1,\x_2\in\atm$, and
\begin{equation*}\tag{$*$}
  \<\u_0,\x_1,\u_2\>,\<\u_2,\x_0,\u_1\>,\<\u_0,\x_2,\u_1\>\in\T'.
\end{equation*}
We wish to show $\x_1\rp\x_0\geq\x_2$. If
$\w\notin\{\u_0,\u_1,\u_2\}$, then $\{\u_0,\u_1,\u_2\}\subseteq\U$, so
we get $\x_1\rp\x_0\geq\x_2$ from the fact that (v) holds in
$\<\U,\T\>$.  So we may assume $\w\in\{\u_0,\u_1,\u_2\}.$ Suppose, for
example, that $\w=\u_0$.  Then ($*$) implies that
$\<\w,\x_1,\u_2\>,\<\w,\x_2,\u_1\>\in\T'.$ Now there are no triples in
$\T$ having $\w$ as first term, so
$\<\w,\x_1,\u_2\>,\<\w,\x_2,\u_1\>\in\T'\minus\T,$ and hence
$\u_1,\u_2\in\{\u,\v,\w\}.$ Arguing similarly in case $\w=\u_1$ or
$\w=\u_2$, we conclude that $\{\u_0,\u_1,\u_2\}\subseteq\{\u,\v,\w\}.$
There are 19 ways of choosing $\u_0,\u_1,\u_2$ so that
$\w\in\{\u_0,\u_1,\u_2\}\subseteq\{\u,\v,\w\}$, and these ways are
listed in the first 3 columns of the table below.

The next 3 columns show which atoms must be assigned to
$\x_0,\x_1,\x_2$ in order that ($*$) holds. These entries are easily
determined from the definition of $\T'$ when $\w$ is involved. For
example, when $\u_0=\u$ and $\u_1=\w$ in line 2, there is only one
atom $\x_2$ such that $\<\u,\x_2,\w\>\in\T'$, namely $\x_2=\b$.  When
two of $\u_0,\u_1,\u_2$ are either both equal to $\u$, or both equal
to $\v$, then an additional argument is required. Consider, for
example, the case $\u=\u_0=\u_1$ (line 13). By ($*$),
$\<\u,\x_2,\u\>=\<\u_0,\x_2,\u_1\>\in\T$. From this and
$\<\u,\a,\v\>\in T$ we have
$\<\u,\x_2,\u\>,\<\u,\a,\v\>,\<\u,\a,\v\>\in\T$, and hence
$\x_2\rp\a\geq\a$ by (v). Thus $\x_2\rp\a\neq0$. We also have
$\x_2\leq\id$ by (iii), and consequently $\x_2=\dom\a$ by
\eqref{domrng} and \eqref{1(39)}.

The seventh column contains the conclusions which must be proved in
each of the 19 cases. Following the table are brief remarks on how
each of the entries in column 7 can be derived.
\begin{equation*} 
\begin{array}{|r|c|c|c|c|c|c|c|}\hline
  &\,\,\,1\,\,\,&\,\,\,2\,\,\,&\,\,\,3\,\,\,&\,\,\,4\,\,\,&\,\,\,5\,\,\,&\,\,\,6\,\,\,&7\\\hline
  &\u_0&\u_1&\u_2&\x_2&\x_1&\x_0&\,\x_2\leq\x_1\rp\x_0\,\\\hline
  1&\u&\v&\w&\a&\b&\c&\a\leq\b\rp\c\\
  2&\u&\w&\v&\b&\a&\con\c&\b\leq\a\rp\con\c\\
  3&\v&\w&\u&\con\c&\con\a&\b&\con\c\leq\con\a\rp\b\\
  4&\v&\u&\w&\con\a&\con\c&\con\b&\con\a\leq\con\c\rp\con\b\\
  5&\w&\u&\v&\con\b&\c&\con\a&\con\b\leq\c\rp\con\a\\
  6&\w&\v&\u&\c&\con\b&\a&\c\leq\con\b\rp\a\\
  7&\u&\w&\w&\b&\b&\rng\b&\b\leq\b\rp\rng\b\\
  8&\w&\u&\w&\con\b&\rng\b&\con\b&\con\b\leq\rng\b\rp\con\b\\
  9&\w&\w&\u&\rng\b&\con\b&\b&\rng\b\leq\con\b\rp b\\
  10&\v&\w&\w&\con\c&\con\c&\rng\b&\con\c\leq\con\c\rp\rng\b\\
  11&\w&\v&\w&\c&\rng\b&\c&\c\leq\rng\b\rp\c\\
  12&\w&\w&\v&\rng\b&\c&\con\c&\rng\b\leq \c\rp\con\c\\
  13&\u&\u&\w&\dom\a&\b&\con\b&\dom\a\leq \b\rp\con\b\\
  14&\u&\w&\u&\b&\dom\a&\b&\b\leq\dom\a\rp\b\\
  15&\w&\u&\u&\con\b&\con\b&\dom\a&\con\b\leq\con\b\rp\dom\a\\
  16&\v&\v&\w&\rng\a&\con\c&\c&\rng\a\leq\con\c\rp c\\
  17&\v&\w&\v&\con\c&\rng\a&\con\c&\con\c\leq\rng\a\rp\con\c\\
  18&\w&\v&\v&\c&\c&\rng\a&\c\leq\c\rp\rng\a\\
  19&\w&\w&\w&\rng\b&\rng\b&\rng\b&\rng\b\leq\rng\b\rp\rng\b\\\hline
\end{array}
\end{equation*}
From the fact that $\<\u,\a,\v,\b,\c\>$ is a flaw in $\<\U,\T\>$ we
get $\a\leq\b\rp\c$, which is the conclusion needed for line 1.  This
is equivalent to $0\neq\b\rp\c\cdot\a$ since $\a$ is an atom.
Applying \eqref{cycle2}, we get $0\neq\a\rp\con\c\cdot\b$, which is
equivalent to $\b\leq\a\rp\con\c$ since $b$ is an atom, thus
accounting for line 2. Each of lines 3 through 6 is obtained in a
similar fashion from the preceding lines. Line 7 follows immediately
from \eqref{zx1867}, line 8 from \eqref{zx1867}, \eqref{ra9},
and \eqref{zx1849}, and line 9 from \eqref{domrng}. Next we show that
lines 10 through 18 can be equivalently transformed into formulas
involving only a single atom.  From lines 1, 2, and 6 we get
$0\neq\b\rp\c$, $0\neq\a\rp\con\c$, and $0\neq\con\b\rp\a$. We
apply \eqref{zx1891} and \eqref{condomrng}, obtaining
$\rng\b=\rng{\con\c}$, $\rng\a=\rng\c$, and $\dom\a=\dom\b$.  Now
substitute $\rng{\con\c}$ for $\rng\b$ in lines 10 to 12, $\dom\b$ for
$\dom\a$ in 13 to 15, and $\rng\c$ for $\rng\a$ in 16 to 18. Then
lines 10 to 18 can be proved like lines 7 to 9.  Finally we note that
line 19 follows from \eqref{zx1850}.

This completes the verification of (v), so the lemma is proved.
\endproof
\begin{lemma}\label{B} Let $\<\U,\T\>$ be a labelling system. Then
  there is a labelling system $\<\U',\T'\>$, which extends
  $\<\U,\T\>$, such that there is no flaw $\<\u,\a,\v,\b,\c\>$ in
  $\<\U',\T'\>$ with $\u,\v\in\U$.
\end{lemma}
\proof Let
$\<\<u_\kappa,a_\kappa,v_\kappa,b_\kappa,c_\kappa\>:\kappa<\alpha\>$
be an enumeration of the flaws in $\<\U,T\>$. Set
$\<\U_0,T_0\>=\<\U,T\>$.  Assume the labelling system
$\<\U_\kappa,T_\kappa\>$ has been constructed.  If
$\<u_\kappa,a_\kappa,v_\kappa,b_\kappa,c_\kappa\>$ is a flaw in
$\<\U_\kappa,T_\kappa\>$, let $\<\U_{\kappa+1},T_{\kappa+1}\>$ be a
labelling system extending $\<\U_\kappa,T_\kappa\>$ in which
$\<u_\kappa,a_\kappa,v_\kappa,b_\kappa,c_\kappa\>$ is not a flaw, by
Lemma \ref{A}, and otherwise let
$\<\U_{\kappa+1},T_{\kappa+1}\>=\<\U_\kappa,T_\kappa\>$.  If
$\lambda\leq\alpha$ is a limit ordinal, then let
\begin{equation*}
  \<\U_\lambda,T_\lambda\>=\Big<\bigcup_{\kappa<\lambda}
    \U_\kappa,\bigcup_{\kappa<\lambda}T_\kappa\Big>.
\end{equation*}
Finally, set
$\<\U',T'\>=\<\U_\alpha,T_\alpha\>$. It is easy to see that $\<\U',T'\>$
is the desired labelling system.
\endproof
\begin{lemma}\label{C}
  Every labelling system $\<\U,T\>$ can be extended to a complete
  labelling system.
\end{lemma}
\proof Set $\<\U_0,T_0\>=\<\U,T\>$, and for each $n<\omega$ ($\omega$
is the least infinite ordinal), let $\<\U_{n+1},T_{n+1}\>$ be the
labelling system obtained by applying Lemma \ref{B} to the previously
constructed labelling system $\<\U_n,T_n\>$.  Then
\begin{equation*}
  \<\U',\T'\>=\Big<\bigcup_{n<\omega}\U_n,\bigcup_{n<\omega}\T_n\Big>
\end{equation*} 
is a complete labelling system extending $\<\U,\T\>$. \endproof

According to the remarks preceding Lemma \ref{A}, the proof of
Theorem~\ref{(10)} is complete.  \endproof

On the beach in 1984, I encouraged Hajnal and Istv\'an to read the
proof just presented, rather than the considerably different method
used in \cite{MR662049} to prove the following corollary.
\begin{cor}[{\rm\cite[Theorem 5.20]{MR662049}}]\label{charWA}
  Let $\K$ be the class of algebras isomorphic to
  $\Rl_\E\big(\Re(\U)\big)$ for some set $\U$ and some symmetric and
  reflexive relation $\E\subseteq\U\times\U$.  Then $\WA$ is the
  closure of $\K$ under the formation of subalgebras:
  $\WA=\mathbf{S}\K.$
\end{cor}
\proof $\K\subseteq\WA$ by Corollary \ref{th1} and
$\WA\subseteq\mathbf{S}\K$ by Theorem \ref{(10)}.  We have
$\mathbf{S}\WA=\WA$ since $\WA$ has equational axiomatizations.
Therefore $\WA\subseteq\mathbf{S}\K\subseteq\mathbf{S}\WA=\WA$.
\endproof
\section{The Relativized Cylindric Set Algebra %
  of a Suitable Structure}\label{sect5}
By Theorem \ref{(10)}, every complete atomic $\gc\A\in\WA$ is a
subalgebra of a relativized set relation algebra.  This result can be
sharpened. By Theorem \ref{th7} in \SS\ref{sect10}, every complete
atomic $\gc\A\in\WA$ is already itself a relativized set relation
algebra, not just a subalgebra of one.  For the proof we will employ
suitable structures and their canonical relativized cylindric set
algebras from \cite{MR987611}.  According to \cite[Definition
1]{MR987611}, $\gc\B=\<\B,\T_\k,\E_{\k\l}\>_{\k,\l<\alpha}$ is a {\bf
  suitable structure} if $\B$ is a set, $\alpha$ is a non-zero
ordinal, and, for all $\k,\l,\m<\alpha$,
\begin{enumerate}
\item[(i)] $\T_\k\subseteq\B\times\B$, $\E_{\k\l}\subseteq\B$,
\item[(ii)] $\T_\k$ is an equivalence relation on $\B$,
\item[(iii)] $\E_{\k\k}=\B$,
\item[(iv)] $\E_{\k\l}=\T_\m^*(\E_{\k\m}\cap\E_{\m\l})$ whenever
  $\k,\l\neq\m$,
\item[(v)] $\T_\k\cap(\E_{\k\l}\times\E_{\k\l})\subseteq Id$ whenever
  $\k\neq\l$.
\end{enumerate}
These are exactly the conditions in \cite[2.7.40]{HMT71}, except the
requirement that $\T_\kappa$ and $\T_\lambda$ commute is deleted.
Assume that $\gc\B=\<\B,\T_\k,\E_{\k\l}\>_{\k,\l<\alpha}$ is a
suitable structure.  Let $\Tr(\gc\B)$ be the set of all sequences
$\p=\<\t_0,\k_0,\dots,\t_\n,\k_\n\>$ such that $\n\in\omega$,
$\t_0,\dots,\t_\n\in\B$, $\k_0,\dots,\k_\n<\alpha$, and, for all
$\i<\n$, $\t_\i\neq\t_{\i+1}$ and $\t_\i\,\T_{\k_i}\,\t_{\i+1}$.  If
$\p\in Tr(\gc\B)$ then $\p$ is called a $\gc\B$-{\bf trail}, $\p$ {\bf
  begins at} $\t_0$, $\p$ {\bf ends at} $\t_\n$, $\k_\n$ is the {\bf
  pointer of} $\p$, $\p$ has {\bf length} $|\p|=\n+1$, and $\p$ is
{\bf reduced} if the following conditions hold:
\begin{enumerate}
\item[(i)] if $1=|\p|$ and $\t_0\in\E_{\k_0\l}$, then
  $\k_0\leq\l<\alpha$,
\item[(ii)] if $1<|\p|$ then $\k_{\n-1}=\k_\n$ and for all
  $\l<\alpha$, $\t_\n\in\E_{\k_n\l}$ iff $\k_\n=\l$,
\item[(iii)] if $0\leq\i<|\p|-2$, then either $\t_\i\neq\t_{\i+2}$ or
  $\k_\i\neq\k_{\i+1}$.
\end{enumerate}
For every $\l<\alpha$, let
$\p\l=\<\t_0,\k_0,\dots,\t_{\n-1},\k_{\n-1},\t_\n,\l\>.$ Let $\approx$
be the smallest equivalence relation on $\Tr(\gc\B)$ such that
\begin{enumerate}
\item[(i)]
  $\<\t_0,\k_0,\dots,\t_i,\l,\s,\l,\t_\i,\k_\i,\dots,\t_\n,\k_\n\>
  \approx\<\t_0,\k_0,\dots,\t_\i,\k_\i,\dots,\t_\n,\k_\n\>$ \newline
  where $0\leq\i\leq\n$,
\item[(ii)] $\<\t_0,\k_0,\dots,\t_\n,\l,\s,\k_\n\>
  \approx\<\t_0,\k_0,\dots,\t_\n,\k_\n\>$ where $\l\neq\k_\n$,
\item[(iii)]
  $\<\t_0,\k_0,\dots,\t_\n,\l\>\approx\<\t_0,\k_0,\dots,\t_\n,\k_\n\>$
  where $\t_\n\in\E_{\l\k_\n}$.
\end{enumerate}
For each $\p\in\Tr(\gc\B)$, let $\p^{\gc\B}$ be the
$\approx$-equivalence class of $\p$, \ie,
$\p^{\gc\B}=\lb\p':\p\,\approx\,\p'\rb.$ The conditions defining
$\approx$ are called ``reductions'', \eg, a reduction of type (i)
consists of the replacement of any subsequence of the form
$\<\t,\l,\s,\l,\t\>$ by $\<\t\>$. Equivalent trails have the same
beginnings, but may have different ends and different pointers, due to
reductions of types (ii) and (iii).  By \cite[Lemma 6]{MR987611},
there is a function called $\P_7$ such that $\P_7(\p)$ is the unique
reduced trail in $\p^{\gc\B}$.  Suppose $\q\in\p^{\gc\B}$ and $\q$ is
reduced. Then $|\q|\leq|\p|$, and if $|\p|=1$, either $\p$ is already
reduced, or else $\q=\p\k$ for some ordinal $\k$ which is strictly
smaller than the pointer of $\p$.  Let
$\U(\gc\B)=\lb\p^{\gc\B}:\p\in\Tr(\gc\B)\rb.$ For every point
$\u\in\U(\gc\B)$, let $|\u|$ be the length of the unique reduced trail
in $\u$.  For every $\gc\B$-trail
$\p=\<\t_0,\k_0,\t_1,\k_1,\dots,\t_{\n-1},\k_{\n-1},\t_\n,\k_\n\>$ let
$\con\p=\<\t_n,\k_{n-1},\t_{n-1},\k_{n-2},\dots,\t_1,\k_0,\t_0,\k_n\>.$
If $\q=\<\s_0,\l_0,\dots,\s_m,\l_m\>$ is any other $\gc\B$-trail,
then $\p\odot\q$ is defined iff $\p$ ends where $\q$ begins, in which
case $\p\odot\q=\<\t_0,\k_0,\t_1,\k_1,\dots,\t_{\n-1},
\k_{\n-1},\s_0,\l_0,\dots,\s_m,\l_m\>.$ Also,
\begin{equation*}
  \L_\p(\q)=\begin{cases}\p\odot\q&\text{ if }\t_\n=\s_0\\
    \con\p\odot\q&\text{ if }\t_\n\neq\s_0=\t_0\\
    \q&\text{ if }\t_\n\neq\s_0\neq\t_0,
  \end{cases}
\end{equation*} 
and, for any $X\subseteq\Tr(\gc\B)$,
$\ell_\p(\X)=\bigcup_{\q\in\X}\(\L_\p(\q)\)^{\gc\B}.$ Let
\begin{equation*}
  Pm(\gc\B)=\lb\ell_\p:\p\in \Tr(\gc\B)\rb.
\end{equation*}
By \cite[Lemma 11]{MR987611}, $Pm(\gc\B)$ is a group of permutations
of $\U(\gc\B)$, and the inverse of the permutation $\ell_\p\in
Pm(\gc\B)$ is $\ell_{\con\p}$.  For every $\t\in\B$, define a set of
sequences of length $\alpha$, by
\begin{equation*}
  \rep\t=\lb\left<(\p\k)^{\gc\B}:\k<\alpha\right>:
  \p\in\Tr(\gc\B),\ \p\text{ ends at }\t\rb,
\end{equation*}
and set $\V(\gc\B)=\bigcup_{\t\in\B}\rep\t.$ By
\cite[3.1.2(i)]{HMT85}, $\Sb({}^\alpha\U(\gc\B))$ is the full
$\alpha$-dimensional cylindric set algebra of all $\alpha$-ary
relations on $\U(\gc\B)$.  If $\t\in\B$ then, by \cite[Lemma
  12(i)]{MR987611}, $Pm(\gc\B)$ preserves $\rep\t$ and acts
transitively on $\rep\t$, hence $\rep\t$ is the orbit, under the
action of the group $Pm(\gc\B)$, of any single $\alpha$-sequence in
$\rep\t$.  It follows, by \cite[Lemma 9]{MR987611}, that $\rep\t$ is
an atom of the subalgebra $\gc\C$ of $\Sb({}^\alpha\U(\gc\B))$
completely generated by $\lb\rep\t:\t\in\B\rb$.  This complete
subalgebra $\gc\C$ is called the {\bf canonical cylindric set algebra}
of the suitable structure $\gc\B$. Its set of atoms includes
$\lb\rep\t:\t\in\B\rb$, but is typically much larger. By relativizing
(in the cylindric algebraic sense of \cite[2.2.1]{HMT71}) to the union
$\V(\gc\B)$ of this set of atoms, one gets
$\Rc\gc\B=\Rl_{\V(\gc\B)}\(\gc\C\)$, the {\bf canonical relativized
  cylindric set algebra} \cite[Definition 7]{MR987611}.  Its set of
atoms is exactly $\lb\rep\t:\t\in\B\rb$ by \cite[Lemma
  12(ii)]{MR987611}.
\section{The Suitable Structure of a $\WA$}\label{sect6}
We show next that every complete atomic weakly associative relation
algebra $\gc\A$ has a suitable structure $\gc\B$ that is
3-dimensional.
\begin{lemma}\label{suitable}
  Given an atomic $\gc\A\in\WA$, let
  $\B=\{\s:\s\in{}^3\atm,\,\s_2\rp\s_0\geq\s_1\}$ and, assuming
  $\{\k,\l,\m\}=\{0,1,2\}$, $\T_\k=\{\<\s,\t\>:\s_\k=\t_\k\}$,
  $\E_{\k\k}=\B$, and
  $\E_{\k\l}=\E_{\l\k}=\{\s:\s\in\B,\,\s_\m\leq\id\}$.  Then
  $\gc\B=\<\B,\T_\k,\E_{\k\l}\>_{\k,\l<3}$ is a suitable structure.
\end{lemma}
\proof The first three properties (i)--(iii) required of suitable
structures clearly hold. We need only show (iv) and (v).  First,
observe that $\B$ is the union of {\bf cycles}, sets of the form
\begin{equation*}
  \C(\b,\c,\a)=\{\<\b,\c,\a\>,\<\c,\b,\con\a\>,\<\con\c,\con\a,\b\>,
  \<\con\a,\con\c,\con\b\>, \<\a,\con\b,\con\c\>,
  \<\con\b,\a,\c\>\},
\end{equation*}
where $\a,\b,\c\in\atm$. To see this, note that if
$\<\b,\c,\a\>\in\B$, then $\con\a,\con\b,\con\c\in\atm$ by
\eqref{1(36)}, and $\a\rp\b\geq\c$, so by \eqref{cyclelaw},
\begin{align*}
  \con\a\rp\c&\geq\b,&\b\rp\con\c&\geq\con\a,&\con\b\rp\con\a
  &\geq\con\c,&\con\c\rp\a&\geq\con\b,&\c\rp\con\b&\geq\a,
\end{align*}
hence $\C\subseteq\B$. Because of axiom \eqref{ra7}, this shows that
$\C\subseteq\B$ whenever any of the listed triples of $\C$ is in $\B$.
$\C$ is a {\bf diversity cycle} if
$\a,\b,\c,\con\a,\con\b,\con\c\leq\di$, and an {\bf identity cycle}
otherwise. 

If $\<\b,\c,\a\>\in\B$ and any two of $\{\a,\b,\c\}$ are below $\id$,
then so is the third, and $\a=\b=\c$. (Their converses are also below
$\id$, so by one of the inequalities above the third element is below
the product of two elements below $\id$, and hence is itself below
$\id$.)

Assume $\<\b,\c,\a\>\in\B$ and one of $\{\a,\b,\c\}$ is below $\id$.
From the inequalities above we can use \eqref{zx1849}, \eqref{zx1890},
\eqref{1(39)}, and definition \eqref{domrng} of $\dom\x$ and $\rng\x$
to deduce that
\begin{align*}\tag{$**$} 
  \<\b,\c,\a\>=\begin{cases}
    \<\b,\b,\dom\b\>=\<\c,\c,\dom\c\>&\text{if }\a\leq\id,\\
    \<\rng\a,\a,\a\>=\<\rng\c,\c,\c\>&\text{if }\b\leq\id,\\
    \<\con\a,\dom\a,\a\>=\<\b,\rng\b,\con\b\> &\text{if }
    \c\leq\id,
\end{cases}
\end{align*}
hence
\begin{align*}
  \C(\b,\c,\a)=\begin{cases}
    \{\<\b,\b,\dom\b\>,\<\con\b,\dom\b,\b\>,\<\dom\b,
    \con\b,\con\b\>\}
    &\text{if }\a\leq\id,\\
    \{\<\rng\a,\a,\a\>,\<\a,\rng\a,\con\a\>,\<\con\a,\con\a,\rng\a\>\}
    &\text{if }\b\leq\id,\\
    \{\<\con\a,\dom\a,\a\>,\<\dom\a,\con\a,\con\a\>,\<\a,\a,\dom\a\>\}
    &\text{if }\c\leq\id.
\end{cases}
\end{align*}
Proof of property (iv). Assume $\t\in\T_\m^*(E_{\k\m}\cap E_{\m\l})$
in the case where $\{\k,\l,\m\}=\{0,1,2\}$.  Then, for some $\s\in\B$,
$\t\mathrel{\T_\m}\s$ and $\s\in\E_{\k\m}\cap\E_{\m\l}$, which gives
us $\t_\m=\s_\m$, $\s_\l\leq\id$, and $\s_\k\leq\id$.  Since two
elements of $\s$ are below $\id$, the third one is as well, \ie,
$\s_\m\leq\id$, hence $\t\in\E_{\k\m}$ since $\t_\m=\s_\m\leq\id$.
This proves the equality in one direction. For the other, assume
$\t\in\E_{\k\m}$, so $\t_\l\leq\id$. We want some
$\s\in\E_{\k\m}\cap\E_{\m\l}$ such that $\s\mathrel{\T_\m}\t$, \ie,
$\t_\m=\s_\m$, $\s_\l\leq\id$, and $\s_\k\leq\id$. It suffices to let
$\s=\<\t_\l,\t_\l,\t_\l\>$.

Since $\E_{\k\m}=\E_{\m\k}$ and $\E_{\k\k}=\B$, the remaining case
of (iv) is $\B=T_\m^*E_{\k\m}$ whenever $\k\neq\m$.  The inclusion
from right to left is trivially true.  Let $\t=\<\b,\c,\a\>\in\B$.
Depending on the values of $\k$ and $\m$, the following list contains
an $\s\in\E_{\k\m}$ such that $\s\mathrel{\T_\m}\t$:
\begin{align*}
&\t=\<\b,\c,\a\>\mathrel{\T_0}\<\b,\b,\dom\b\>\in\E_{01},&&
&\t=\<\b,\c,\a\>\mathrel{\T_0}\<\b,\rng\b,\con\b\>\in\E_{02},\\
&\t=\<\b,\c,\a\>\mathrel{\T_1}\<\c,\c,\dom\c\>\in\E_{01},&&
&\t=\<\b,\c,\a\>\mathrel{\T_1}\<\dom\c,\c,\c\>\in\E_{12},\\
&\t=\<\b,\c,\a\>\mathrel{\T_2}\<\con\a,\dom\a,\a\>\in\E_{02},&&
&\t=\<\b,\c,\a\>\mathrel{\T_2}\<\rng\a,\a,\a\>\in\E_{12}.
\end{align*}
Proof of property (v). Assume
$\<\s,\t\>\in\T_\k\cap(\E_{\k\l}\times\E_{\k\l})$ where $\k\neq\l$.
We wish to show $\s=\t$.  From the hypothesis we get
$\s\mathrel{\T_\k}\t$ and $\s,\t\in\E_{\k\l}$, \ie, $\s_\k=\t_\k$,
$\s_\m\leq\id$, and $\t_\m\leq\id$, where $\{0,1,2\}=\{\k,\l,\m\}$.
Since $\s$ and $\t$ are triples with a subidentity element in the same
position $\m$, they must have the same form according to ($**$), and
since they also have the same element at a different position $\k$,
they are the same.  For example, if $\m=2$, then
$\s=\<\s_0,\s_0,\dom{\s_0}\>$ and $\t=\<\t_0,\t_0,\dom{\t_0}\>$. But
$\s_0=\t_0$ follows from $\s\mathrel{\T_\k}\t$ if $\k$ is either 0 or
1, so $\s=\t$.  If $\m=1$ then $\s=\<\breve{\s_2},\dom{\s_2},{\s_2}\>$
and $\t=\<\con{\t_2},\dom{\t_2},{\t_2}\>$.  We get $\s_2=\t_2$ from
$\s\mathrel{\T_2}\t$ and $\breve{\s_2}=\con{\t_2}$ from
$\s\mathrel{\T_0}\t$, and again $\s=\t$.  If $\m=0$ then
$\s=\<\rng{\s_2},\s_2,\s_2\>$ and $\t=\<\rng{\t_2},\t_2,\t_2\>$, so
$\s_2=\t_2$ follows from $\s\mathrel{\T_\k}\t$ if $\k$ is either 1 or
2, and again $\s=\t$.  \endproof
\section{The Complex Algebra of a Suitable Structure}\label{sect7}
If $\gc\B=\<\B,\T_\k,\E_{\k\l}\>_{\k,\l<3}$ is a suitable structure,
then the {\bf complex algebra of $\gc\B$} \cite[2.7.33]{HMT71} is
\begin{equation*}
  \Cm\gc\B=\{\pow(\B),\cup,\cap,{}_\B{\minus},\emptyset,\B,
  \T^*_\i,\E_{\i\j}\}_{\i,\j<3},
\end{equation*}
where, for every $\X\subseteq\B$, $\T^*_\i(\X)=\{\y:\y\T_\i\x\in\X\}$.
The class $\NCA_\alpha$ is defined by the axioms for
$\alpha$-dimensional cylindric algebras $\CA_\alpha$ with the
commutativity of cylindrifications deleted, that is, by the axioms
(C$_0$)--(C$_3$), (C$_5$)--(C$_7$) of \cite[1.1.1]{HMT71}.  The class
$\NA_\alpha$ of {\bf non-commutative cylindric algebras of dimension
  $\alpha$} is defined by the axioms of $\CA_\alpha$ with postulate
(C$_4$) of \cite[1.1.1]{HMT71}, $\cyl_\k\cyl_\l x =\cyl_\l\cyl_\k x$,
replaced by the weaker postulate (C$_4^*$), $\cyl_\k\cyl_\l
x\geq\cyl_\l\cyl_\k x\cdot\diag_{\l\m}$ with $\m\neq\k,\l$ (see
\cite{AT}, \cite{N86}, and \cite{T}). The next lemma is a sharpening
of one direction of \cite[Lemma 15]{MR662049}.  It is about $\NA_3$
instead of the wider class $\NCA_3$.
\begin{lemma}\label{NA3}
  Assume $\gc\A\in\WA$, $\gc\A$ is atomic, and
  $\gc\B=\<\B,\T_\k,\E_{\k\l}\>_{\k,\l<3}$ is the suitable structure
  built from $\gc\A$ in Lemma \ref{suitable}.  Then
  $\Cm\gc\B\in\NA_3$.
\end{lemma}
\proof By \cite[Lemma 15]{MR987611}, $\Cm\gc\B\in\NCA_3$, so we need
only show $\Cm\gc\B$ satisfies (C$_4^*$). In the notation of complex
algebras, (C$_4^*$) says that if $\X\subseteq\B$ and $\m\neq\k,\l$,
then $\T^*_\l\T^*_\k\X\cap\E_{\l\m}\subseteq\T^*_\k\T^*_\l\X.$ Since
this is trivially true if $\k=\l$, we assume $\k\neq\l$, hence
$\{0,1,2\}=\{\k,\l,\m\}$.  Suppose $\t\in\B$ and
$\t=\<\b,\c,\a\>\in\T^*_\l\T^*_\k\X\cap\E_{\l\m}.$ Then $\t_\k\leq\id$
since $\t\in\E_{\l\m},$ so, by ($**$),
\begin{align*}
  \t=\begin{cases}
    \<\b,\b,\dom\b\>&\text{if }\k=2,\\
    \<\rng\a,\a,\a\>&\text{if }\k=0,\\
    \<\con\a,\dom\a,\a\>&\text{if }\k=1.
\end{cases}
\end{align*}
From $\t\in\T^*_\l\T^*_\k\X$ we know there are $\s,\r\in\B$ such that
$\t\mathrel{\T_\l}\s\mathrel{\T_\k}\r\in\X.$ We wish to find, in each
of the six cases, some $\s'\in\B$ such that
$\t\mathrel{\T_\k}\s'\mathrel{\T_\l}\r\in\X.$ The following table
shows how to compute $\s'$ from the values of $\k,\l,\m$. It is
immediately obvious from the definition of $\s'$ that
$\s'\mathrel{\T_\l}\r$, while $\t\mathrel{\T_\l}\s'$ follows from the
reason given in the last column, deduced for each case right after the
table.
\begin{align*}
  \<\k,&\l,\m\>=&&\t=&&\s'=
  &\t\mathrel{\T_\l}\s'&\text{ because}\\
  &\<2,0,1\>&&\<\b,\b,\dom\b\> &&\<\r_0,\r_0,\dom{\r_0}\>
  &&\dom\b=\dom{\r_0} \\
  &\<2,1,0\>&&\<\b,\b,\dom\b\> &&\<\r_1,\r_1,\dom{\r_1}\>
  &&\dom\b=\dom{\r_1} \\
  &\<0,1,2\>&&\<\rng\a,\a,\a\> &&\<\rng{\r_1},\r_1,\r_1\>
  &&\rng\a=\rng{\r_1} \\
  &\<0,2,1\>&&\<\rng\a,\a,\a\> &&\<\rng{\r_2},\r_2,\r_2\>
  &&\rng\a=\rng{\r_2} \\
  &\<1,0,2\>&&\<\con\a,\dom\a,\a\>&&\<\r_0,\rng{\r_0},\con{\r_0}\>
  &&\dom\a=\rng{\r_0} \\
  &\<1,2,0\>&&\<\con\a,\dom\a,\a\>&&\<\con{\r_2},\dom{\r_2},\r_2\>
  &&\dom\a=\dom{\r_2}
\end{align*}
From $\s,\r\in\B$ it follows, by \eqref{cyclelaw}, \eqref{condomrng},
and \eqref{zx1891}, that statements (a)--(f) hold:
(a) $\dom{\s_1}=\dom{\s_2}$, (b)\ $\rng{\s_2}=\dom{\s_0}$, (c)\
$\rng{\s_0}=\rng{\s_1}$, (d)\ $\dom{\r_1}=\dom{\r_2}$, (e)\
$\rng{\r_2}=\dom{\r_0}$, and (f)\ $\rng{\r_0}=\rng{\r_1}$.  Each
equation in the last column of the table above may be deduced as
follows, using various combinations of statements (a)--(f) with the
hypotheses $\t\mathrel{\T_\l}\s\mathrel{\T_\k}\r$.
\begin{align*}
  \<2,0,1\>&&\dom\b=\dom{\t_0}\overset{\t\mathrel{\T_0}\s}=
  \dom{\s_0}&\overset{\text{(b)}}=\rng{\s_2}\overset{\s\mathrel{\T_2}\r}=
  \rng{\r_2}\overset{\text{(e)}}=\dom{\r_0},\\
  \<2,1,0\>&&\dom\b=\dom{\t_1}\overset{\t\mathrel{\T_1}\s}=
  \dom{\s_1}&\overset{\text{(a)}}=\dom{\s_2}\overset{\s\mathrel{\T_2}\r}=\dom{\r_2}
  \overset{\text{(d)}}=\dom{\r_1},\\
  \<0,1,2\>&&\rng\a=\rng{\t_1}\overset{\t\mathrel{\T_1}\s}=
  \rng{\s_1}&\overset{\text{(c)}}=\rng{\s_0}\overset{\s\mathrel{\T_0}\r}=\rng{\r_0}
  \overset{\text{(f)}}=\rng{\r_1},\\
  \<0,2,1\>&&\rng\a=\rng{\t_2}\overset{\t\mathrel{\T_2}\s}=
  \rng{\s_2}&\overset{\text{(b)}}=\dom{\s_0}\overset{\s\mathrel{\T_0}\r}=\dom{\r_0}
  \overset{\text{(e)}}=\rng{\r_2},\\
  \<1,0,2\>&&\dom\a\overset{\eqref{condomrng}}=\,\rng{\con\a}=\rng{\t_0}
  \overset{\t\mathrel{\T_0}\s}=\rng{\s_0}&\overset{\text{(c)}}=\rng{\s_1}
  \overset{\s\mathrel{\T_1}\r}=\rng{\r_1}\overset{\text{(f)}}=\rng{\r_0},\\
  \<1,2,0\>&&\dom\a=\dom{\t_2}\overset{\t\mathrel{\T_2}\s}=
  \dom{\s_2}&\overset{\text{(a)}}=\dom{\s_1}\overset{\s\mathrel{\T_1}\r}=
  \dom{\r_1}\overset{\text{(d)}}=\dom{\r_2}.
\end{align*}
\endproof
Let $2\leq\n<\omega$.  MGR$_\n$ \cite[3.2.88(1)]{HMT85} is the set of
{\bf $\n$-ary merry-go-round identities}, where
$\l,\k_1,\dots,\k_\n<\alpha$ are distinct ordinals:
\begin{equation*}
  \sub^\l_{\k_1}\sub^{\k_1}_{\k_2}\sub^{\k_2}_{\k_3}\dots
  \sub^{\k_{\n-1}}_{\k_\n}\sub^{\k_\n}_{\l}\cyl_\l\x=
  \sub^\l_{\k_\n}\sub^{\k_\n}_{\k_1}\sub^{\k_1}_{\k_2}
  \dots\sub^{\k_{\n-2}}_{\k_{\n-1}}\sub^{\k_{\n-1}}_{\l}\cyl_\l\x.
\end{equation*}
In particular, MGR$_2$ \cite[3.2.88(2)]{HMT85} is the set of $2$-ary
merry-go-round identities: $$
\sub^\l_{\k_1}\sub^{\k_1}_{\k_2}\sub^{\k_2}_{\l}\cyl_\l x =
\sub^\l_{\k_2}\sub^{\k_2}_{\k_1}\sub^{\k_1}_{\l}\cyl_\l x,$$ and
MGR$_3$ \cite[3.2.88(3)]{HMT85} is the set of $3$-ary merry-go-round
identities: $$\sub^\l_{\k_1}\sub^{\k_1}_{\k_2}\sub^{\k_2}_{\k_3}
\sub^{\k_3}_{\l}\cyl_\l\x=
\sub^\l_{\k_3}\sub^{\k_3}_{\k_1}\sub^{\k_1}_{\k_2}
\sub^{\k_{2}}_{\l}\cyl_\l\x.$$
\begin{lemma}\label{functions} If $\gc\B$ is the suitable structure
  of an atomic $\gc\A\in\WA$, and $\t\in\B$, then
  \begin{enumerate}
  \item[\rm(i)]
    $\E_{10}\cap\T^*_0\{\t\}=\{\<\t_0,\t_0,\dom{\t_0}\>\}$,
  \item[\rm(ii)]
    $\E_{01}\cap\T^*_1\{\t\}=\{\<\t_1,\t_1,\dom{\t_1}\>\}$,
  \item[\rm(iii)]
    $\E_{12}\cap\T^*_2\{\t\}=\{\<\rng{\t_2},\t_2,\t_2\>\}$,
  \item[\rm(iv)]
    $\E_{21}\cap\T^*_1\{\t\}=\{\<\rng{\t_1},\t_1,\t_1\>\}$,
  \item[\rm(v)]
    $\E_{02}\cap\T^*_2\{\t\}=\{\<\con{\t_2},\dom{\t_2},\t_2\>\}$,
  \item[\rm(vi)]
    $\E_{20}\cap\T^*_0\{\t\}=\{\<\t_0,\rng{\t_0},\con{\t_0}\>\}$.
\end{enumerate}
\end{lemma}
\proof The inclusions from right to left follow immediately from
definitions.  The opposite inclusions have very similar proofs, so we
do only the first one.  Assume $\s\in\E_{10}\cap\T^*_0\{\t\}$.  Then
$\s_2\leq\id$ since $\s\in\E_{10}$, so $\s=\<\s_0,\s_0,\dom{\s_0}\>$
by ($**$).  We also have $\s\mathrel{\T_0}\t$ since
$\s\in\T^*_0\{\t\}$, hence $\s_0=\t_0$, yielding
$\s=\<\s_0,\s_0,\dom{\s_0}\> =\<\t_0,\t_0,\dom{\t_0}\>$, as desired.
\endproof
\begin{lemma}\label{MGR} If $\gc\B$ is the suitable structure
  of an atomic $\gc\A\in\WA$ then $\Cm\gc\B$ satisfies {\rm MGR}$_\n$
  for all $\n$.
\end{lemma}
\proof MGR$_\n$ is satisfied in all 3-dimensional algebras whenever
$3\leq\n$, because the index ordinals in MGR$_\n$ are required to be
distinct and dimension 3 is not large enough to find four (or more)
distinct ordinals.  Since $\Cm\gc\B\in\NA_3$ by Lemma \ref{NA3}, we
need only prove MGR$_2$.  In the notation of complex algebras, MGR$_2$
says that for every $\X\subseteq\B$, if $\{\k,\l,\m\}=\{0,1,2\}$ then
$\T^*_\k(\E_{\k\l}\cap\T^*_\l(\E_{\l\m}\cap\T^*_\m(\E_{\m\k}\cap\T^*_\k\X)))
=\T^*_\k(\E_{\k\m}\cap\T^*_\m(\E_{\m\l}\cap\T^*_\l(\E_{\l\k}\cap\T^*_\k\X))).$
Here we compute both sides of a single instance of MGR$_2$, showing
that they evaluate to the same thing.  (The first computation will
also be used in the proof of Theorem \ref{th5}.) The remaining cases
can be obtained from this one by permuting subscripts and making other
appropriate changes. First we evaluate the MGR$_2$ term that comes
from converse.
\begin{align*}
  \con{\{\t\}}
  &=\T^*_2(\E_{20}\cap\T^*_0(\E_{01}\cap\T^*_1(\E_{12}\cap\T^*_2\{\t\})))\\
  &=\T^*_2(\E_{20}\cap\T^*_0(\E_{01}\cap\T^*_1\{\<\rng{\t_2},\t_2,\t_2\>\}))
  &&\text{Lemma \ref{functions}(iii)}\\
  &=\T^*_2(\E_{20}\cap\T^*_0\{\<\t_2,\t_2,\dom{\t_2}\>\})
  &&\text{Lemma \ref{functions}(ii)}\\
  &=\T^*_2\{\<\t_2,\rng{\t_2},\con{\t_2}\>\}
  &&\text{Lemma \ref{functions}(vi)}\\
  &=\{\s:\s\in\B,\,\s_2=\con{\t_2}\}.
\end{align*}
Next we compute the other side of the MGR$_2$ identity connected with
converse. This term can serve as an alternative definition of
converse.
\begin{align*}
  &\T^*_2(\E_{21}\cap\T^*_1(\E_{10}\cap\T^*_0(\E_{02}\cap\T^*_2\{\t\})))\\
  &=\T^*_2(\E_{21}\cap\T^*_1(\E_{10}\cap\T^*_0
  \{\<\con{\t_2},\dom{\t_2},\t_2\>\}))
  &&\text{Lemma \ref{functions}(v)}\\
  &=\T^*_2(\E_{21}\cap\T^*_1\{\<\con{\t_2},\con{\t_2},\dom{\con{\t_2}}\>\})
  &&\text{Lemma \ref{functions}(i)}\\
  &=\T^*_2\{\<\rng{\con{\t_2}},\con{\t_2},\con{\t_2}\>\}
  &&\text{Lemma \ref{functions}(iv)}\\
  &=\{\s:\s\in\B,\,\s_2=\con{\t_2}\}.
\end{align*}
\endproof

\begin{theorem}\label{th4}
  Suppose $\gc\B$ is the suitable structure of an atomic
  $\gc\A\in\WA$. Then $\Cm\gc\B\cong\Rc\gc\B,$ via the isomorphism
  $\rep{}$ from $\Cm\gc\B$ to $\Rc\gc\B$ defined for $\X\subseteq\B$
  by
\begin{align*}
  \rep{}(\X)&=\bigcup_{\t\in\X}\rep\t
  =\bigcup_{\t\in\X}\lb\<(\p\k)^{\gc\B}:\k<3\>:\p\text{ ends at }\t,\,
  \p\in\Tr(\gc\B)\rb.
\end{align*}
\end{theorem}
\proof \cite[Theorem C]{MR987611} says that if $2\leq\alpha$,
$\gc\C\in\NA_\alpha$, $\gc\C$ is complete, atomic, and satisfies
MGR$_2$ and MGR$_3$, then $\gc\C\cong\Rc\mathfrak{At}\gc\C$, where
$\mathfrak{At}\gc\C$ is the {\bf atom structure} of $\gc\C$
\cite[2.7.32]{HMT71}.  This theorem applies to $\gc\C=\Cm\gc\B$,
because $\Cm\gc\B$ is complete, atomic, in $\NA_3$ by Lemma \ref{NA3},
and satisfies the MGR identities by Lemma \ref{MGR}, so
$\Cm\gc\B\cong\Rc\mathfrak{At}\Cm\gc\B$ via the isomorphism
of \cite[Definition 16]{MR987611}. But
$\mathfrak{At}\Cm\gc\B\cong\gc\B$ by \cite[2.7.35]{HMT71}, via the
isomorphism that sends $\{\x\}$ to $\x$. Therefore
$\Cm\gc\B\cong\Rc\gc\B$ via the formula above for the composed
isomorphisms.
\endproof
\section{Relation-Algebraic Reducts}\label{sect8}
Consider an arbitrary non-commutative 3-dimensional cylindric algebra
\begin{equation*}
  \gc\C=\<\C,+,\min\blank,\cyl_\kappa,\diag_{\kappa\lambda}
  \>_{\kappa,\lambda<\alpha}\in\NA_3.
\end{equation*}
Following \cite[2.6.28]{HMT71}, restricted to the case $\alpha=3$, let
\begin{equation*}
  Nr_2\gc\C=\{\x:\x\in\C\text{ and }\cyl_2\x=\x\}.
\end{equation*} 
$Nr_2\gc\C$ is the set of {\bf $2$-dimensional elements} of $\gc\C$.
By \cite[5.3.7]{HMT85} (generalized from $\CA_\alpha$ to
$\NA_\alpha$), the binary operation $\rp$ is defined for all
$\x,\y\in\C$ by
\begin{equation*}
  \x\rp\y=\cyl_2(\cyl_1(\diag_{12}\bp\x)\bp\cyl_0(\diag_{02}\bp\y)),
\end{equation*} 
and the unary operation $\con\blank$ is defined for every $\x\in\C$ by
\begin{equation*}
  \con\x=\cyl_2(\diag_{20}\bp\cyl_0(\diag_{01}\bp\cyl_1(\diag_{12}\bp\x))).
\end{equation*}
We have $\cyl_2\cyl_2\x=\cyl_2\x$ for all $\x\in\C$
by \cite[1.2.3]{HMT71}. This implies the set of 2-dimensional elements
is closed under the operations $\rp$ and $\con\blank$. $Nr_2\gc\C$ is
closed under $+$ because $\cyl_2$ distributes over $+$
by \cite[1.2.6]{HMT71}, and $\diag_{01}$ is a 2-dimensional element
because $\cyl_2\diag_{01}=\diag_{01}$ by \cite[1.3.3]{HMT71}.  The
complement $\min\x$ of $\x\in Nr_2\gc\C$ is also 2-dimensional
by \cite[1.2.12]{HMT71}. We therefore have the
algebra \cite[5.3.7]{HMT85}
\begin{equation*}
  \Ra\gc\C=\<Nr_2\gc\C,\,+,\,\min\blank,\,\rp,\,\con\blank,\,\diag_{01}\>.
\end{equation*}
\begin{theorem}\label{th5}
  If $\gc\B=\<\B,\T_\k,\E_{\k\l}\>_{\k,\l<3}$ is the suitable
  structure of a complete atomic $\gc\A\in\WA$, then
  $\gc\A\cong\Ra\Cm\gc\B$ via the isomorphism $\varphi:\A\to\pow(\B)$
  defined by $\varphi(\x)=\{\t:\t\in\B,\,\t_2\leq\x\}$ for all
  $\x\in\A$.
\end{theorem}
\proof Note that $\T^*_2\varphi(\x)\subseteq\varphi(\x)$ since $\T_2$
is an equivalence relation, and if $\s\in\T^*_2\varphi(\x)$ then, by
the definition of $\varphi$, there is some $\t\in\B$ such that
$\s_2=\t_2$ and $\t_2\leq\x$, hence $\s_2\leq\x$, \ie,
$\s\in\varphi(\x)$. This proves that $\T^*_2\varphi(\x)=\varphi(\x)$
for every $\x\in\A$, so, in fact, $\varphi:\A\to Nr_2\Cm\gc\B,$ and
$\varphi$ is a Boolean homomorphism because
\begin{align*}
  \varphi(1)&=\{\t:\t\in\B,\,\t_2\leq1\}=\B,\\
  \varphi(0)&=\{\t:\t\in\B,\,\t_2\leq0\}=\emptyset,\\
  \varphi(\x+\y)&=\{\t:\t\in\B,\,\t_2\leq\x+\y\}\\
  &=\{\t:\t\in\B,\,\t_2\leq\x\text{ or }\t_2\leq\y\}\\
  &=\{\t:\t\in\B,\,\t_2\leq\x\}\cup\{\t:\t\in\B,\,\t_2\leq\y\}\\
  &=\varphi(\x)\cup\varphi(\y),\\
  \varphi(\min\x)&=\{\t:\t\in\B,\,\t_2\leq\min\x\}\\
  &=\{\t:\t\in\B,\,\t_2\not\leq\x\}\\
  &={}_\B{\minus}\varphi(\x).
\end{align*}
Converse was handled earlier---the first computation in the proof of
Lemma \ref{MGR} happens to also show $\con{\{\t\}}=
\varphi(\con{\t_2})$. If $\X,\Y\in Nr_2\Cm\gc\B$ then
\begin{align*}
  &\T^*_1(\E_{12}\cap\X)\cap\T^*_0(\E_{02}\cap\Y)\\
  &=\T^*_1(\E_{12}\cap\T^*_2\X)\cap\T^*_0(\E_{02}\cap\T^*_2\Y)\\
  &=\T^*_1\{\<\rng{\r_2},\,\r_2,\,\r_2\>:\r\in\X\}
  \cap\T^*_0\{\<\con{\q_2},\,\dom{\q_2},\,\q_2\>:\q\in\Y\}
  &&\text{Lemma \ref{functions}(iii)}\\
  &=\{\<\con{\q_2},\,\r_2,\,\a\>:\r\in\X,\,\q\in\Y,
  \,\<\con{\q_2},\,\r_2,\,\a\>\in\B\}\\
  &=\{\<\con{\q_2},\,\r_2,\,\a\>:\r\in\X,\,\q\in\Y,
  \,\a\rp\con{\q_2}\geq\r_2\}\\
  &=\{\<\con{\q_2},\,\r_2,\,\a\>:\r\in\X,\,\q\in\Y,
  \,\r_2\rp\q_2\geq\a\}&&\eqref{cyclelaw}
\end{align*}so
\begin{align*}
  \varphi(\x)\rp\varphi(\y)&=\T^*_2(\T^*_1(\E_{12}
  \cap\varphi(\x))\cap\T^*_0(\E_{02}\cap\varphi(\y)))\\
  &=\T^*_2(\{\<\con{\q_2},\,\r_2,\,\a\>:\r_2\rp\q_2\geq\a,\,
  \r\in\varphi(\x),\,\q\in\varphi(\y)\})\\
  &=\{\t:\t\in\B,\,\r_2\rp\q_2\geq\t_2,\,\r\in\varphi(\x),\,
  \q\in\varphi(\y)\}\\
  &=\{\t:\t\in\B,\,\r_2\rp\q_2\geq\t_2,\,\r_2\leq\x,\,
  \q_2\leq\y,\,\r,\q\in\B\}\\
  &\subseteq\{\t:\t\in\B,\,\x\rp\y\geq\t_2\}=\varphi(\x\rp\y).
\end{align*}
For the inclusion in the other direction, suppose
$\t\in\varphi(\x\rp\y)$, \ie, $\t\in\B$ and $\x\rp\y\geq\t_2$. By
complete additivity \eqref{ca} there are atoms $\a,\b\in\atm$ such
that $\t_2\leq\a\rp\b$, $\a\leq\x$, and $\b\leq\y$.  Let
$\r=\<\rng\a,\a,\a\>$ and $\q=\<\rng\b,\b,\b\>$ and note that
$\r,\q\in\B$. Then $\r_2=\a\leq\x$ and $\q_2\leq\b\leq\y$, so
$\t\in\varphi(\x)\rp\varphi(\y)$ by the first part of the computation.
\endproof
\section{Cylindric-Relativized Representation}\label{sect9}
\begin{theorem}\label{th6}
  Assume $\gc\A$ is a complete atomic $\WA$, $\gc\B$ is the suitable
  structure of $\gc\A$, and for every $\t\in\B$,
  \begin{equation*}
    \rep\t=\lb\<(\p0)^{\gc\B},(\p1)^{\gc\B},(\p2)^{\gc\B}\>
    :\p\in\Tr(\gc\B),\ \p\text{ ends at }\t\rb.
  \end{equation*}
  Assume $\gc\C$ is the subalgebra of the full 3-dimensional cylindric
  set algebra $\Sb\(\un\)$ completely generated by
  $\lb\rep\t:\t\in\B\rb$.  Then
  \begin{equation*}
    \gc\A\cong\Ra\Rl_{\V(\gc\B)}\(\gc\C\)
  \end{equation*}
  and $\lb\rep\t:\t\in\B\rb$ is the set of atoms of
  $\Rl_{\V(\gc\B)}\(\gc\C\)$.
\end{theorem}
\proof Recall from \SS\ref{sect5}, $\lb\rep\t:\t\in\B\rb$ is the set
of atoms of $\Rl_{\V(\gc\B)}\(\gc\C\)$ by \cite[Lemma
12(ii)]{MR987611}.  Relativizing $\gc\C$ to
$\V(\gc\B)=\bigcup_{\t\in\B}\rep\t$ gives the canonical relativized
cylindric set algebra $\Rc\gc\B$,
\begin{align*}
  \Rc\gc\B&=\Rl_{\V(\gc\B)}\(\gc\C\) =
  \Big\<Rl_{\V(\gc\B)}\(\gc\C\),\cup,{}_{\V(\gc\B)}\setminus,
  \Cyl_\k^{[\V(\gc\B)]},\Diag_{\k\l}^{[\V(\gc\B)]}\Big\>_{\k,\l<3},
  \intertext{where}
  Rl_{\V(\gc\B)}\(\gc\C\)&=\{\X:\V(\gc\B)\supseteq\X\in\C\},\\
  \Cyl_\k^{[\V(\gc\B)]}(\X)&=\{\v\in\V(\gc\B):\v_\k=\u_\k,\,\u\in\X\},\\
  \Diag_{\k\l}^{[\V(\gc\B)]}&=\{\v\in\V(\gc\B):\v_\k=\v_\l\}.
\end{align*}
We have $\Cm\gc\B\cong\Rc\gc\B$ by Theorem \ref{th4}, and
$\gc\A\cong\Ra\Cm\gc\B$ by Theorem \ref{th5}, so $
\gc\A\cong\Ra\Rc\gc\B=\Ra\Rl_{\V(\gc\B)}\(\gc\C\) $ via the isomorphism 
\begin{equation*}
\x\mapsto\lb\<(\p\k)^{\gc\B}:\k<3\>:\p\in\Tr(\gc\B),\
\p\text{ ends at }\t\in\B,\,\t_2\leq\x\rb.
\end{equation*}
\endproof
\section{Relativized Relational Representation}\label{sect10}
Suppose $\gc\B$ is the suitable structure of a complete atomic
$\gc\A\in\WA$.  Given a triple $\<\u,\v,\w\>\in\un$, we say that the
ordered pairs in $\{\u,\v,\w\}\times\{\u,\v,\w\}$ {\bf occur in}
$\<\u,\v,\w\>$, and that $\<\u,\v,\w\>$ {\bf carries} those pairs.
Let $\V_2$ be the set of ordered pairs that occur in triples in
$\V(\gc\B)$, \ie,
$\V_2=\{\<\u,\v\>:\u=(\p\k)^{\gc\B},\,\v=(\p\l)^{\gc\B},\,
\p\in\Tr\gc\B,\,\k,\l<3\}.$ Given a trail $\p\in\Tr\gc\B$, {\bf the
  triple of} $\p$ is $\<(\p0)^{\gc\B},(\p1)^{\gc\B},(\p2)^{\gc\B}\>$.
For every trail $\p\in\Tr(\gc\B)$ we will define a function
$\g_\p:\{\u,\v,\w\}\times\{\u,\v,\w\}\to\atm$, where
$\<\u,\v,\w\>=\<(\p0)^{\gc\B},(\p1)^{\gc\B},(\p2)^{\gc\B}\>$.  Let
$\t\in\B$ be the end of $\p$. Then $\t_2\rp\t_0\geq\t_1\neq0$ since
$\t\in\B$, so $\rng{\t_2}=\dom{\t_0}$ by \eqref{zx1891}.  From
$\t_2\rp\t_0\geq\t_1\neq0$ we get
$\t_0\rp\con{\t_1}\geq\con{\t_2}\neq0$ and
$\con{\t_2}\rp\t_1\geq\t_0\neq0$ by  \eqref{cyclelaw}, so
$\rng{\t_0}=\dom{\con{\t_1}}=\rng{\t_1}$ and
$\dom{\t_2}=\rng{\con{\t_2}}=\dom{\t_1}$ by \eqref{zx1891}
and \eqref{condomrng}.  Then $\g_\p$ is defined by
\begin{align*}
  \g_\p(\u,\u)&=\dom{\t_1}=\dom{\t_2},&\g_\p(\v,\v)&=\rng{\t_2}=\dom{\t_0},
  &\g_\p(\w,\w)&=\rng{\t_0}=\rng{\t_1},\\
  \g_\p(\u,\v)&=\t_2,&\g_\p(\v,\w)&=\t_0,&\g_\p(\u,\w)&=\t_1,\\
  \g_\p(\v,\u)&=\con{\t_2},&\g_\p(\w,\v)&=\con{\t_0},&\g_\p(\w,\u)
  &=\con{\t_1}.
\end{align*}
As a general rule, the end of $\p$ is
$\big\<\g_\p(\v,\w),\g_\p(\u,\w),\g_\p(\u,\v)\big\>$, where
$\<\u,\v,\w\>$ is the triple of $\p$.  Further applications of
\eqref{cyclelaw}, \eqref{condomrng}, \eqref{zx1867}, and
\eqref{zx1891} show that for all $\x,\y,\z\in\{\u,\v,\w\}$, $
\g_\p(\x,\y)\rp\g_\p(\y,\z)\geq\g_\p(\x,\z).$ Since $\g_\p$ depends
only on the end of $\p$, and type (i) reductions do not change the end
or the triple of a trail, it follows that $\g_\p=\g_\q$ for any two
trails $\p\approx\q$ that are related by condition (i).  Suppose two
trails $\p\approx\q$ are related by condition (iii). According to the
condition, they must agree except possibly at the pointer, so
$\p\k=\q\k$ for all $\k<3$, hence $\g_\p=\g_\q$.  Finally, suppose two
trails $\p\approx\q$ are related by a pointer-preserving reduction of
type (ii), \eg,
\begin{align*}
  \p0&=\<\dots,\s,1,\t,0\>\approx\<\dots,\s,0\>=\q0,\\
  \p1&=\<\dots,\s,1,\t,1\>\approx\<\dots,\s,1,\t,1\>,\\
  \p2&=\<\dots,\s,1,\t,2\>\approx\<\dots,\s,2\>=\q2.
\end{align*}
Note that $s_1=\t_1$ since $\p$ is a trail. Adding this equation to
the parts of the definitions of $\g_\q$ and $\g_\p$ that involve index
$1$, we get
\begin{align*}
  \g_\q(\u,\u)&=\dom{\s_1}=\dom{\t_1}=\g_\p(\u,\u),\\
  \g_\q(\w,\w)&=\rng{\s_1}=\rng{\t_1}=\g_\p(\w,\w),\\
  \g_\q(\u,\w)&=\s_1=\t_1=\g_\q(\u,\w),\\
  \g_\q(\w,\u)&=\breve{\s_1}=\con{\t_1}=\g_\q(\w,\u),
\end{align*}
so $\g_\p$ and $\g_\q$ agree where they are both defined.  It follows
by induction on the generation of $\approx$ from reductions of types
(i)--(iii), that if $\p\approx\q$ then $\g_\p$ and $\g_\q$ agree
whenever they are both defined.  Every pair $\<\u,\v\>$ in $\V_2$
occurs in the triple of some trail $\p$, and is thereby assigned to
the atom $\g_\p(\u,\v)$.  If $\<\u,\v\>$ also occurs in the triple of
some trail $\q$, and is thereby assigned to the atom $\g_\q(\u,\v)$,
we have $\u=(\p\k)^{\gc\B}=(\q\k')^{\gc\B}$ and
$\v=(\p\l)^{\gc\B}=(\q\l')^{\gc\B}$, hence $\p\k\approx\q\k'$ and
$\p\l\approx\q\l'$. Since the pointer of $\p$ is irrelevant to the
definition of $\g_\p$, we have, as a rule, $\g_\p=\g_{\p\k}$ for all
$\k<3,$ so from either of these equivalences we conclude that
$\g_\p(\u,\v)=\g_\q(\u,\v)$.  Thus there is a function
$\g:\V_2\to\atm$ such that $\g(\u,\v)=\g_\p(\u,\v)$ for any trail $\p$
whose triple contains $\<\u,\v\>$, with the key property
\begin{align*}\tag{K}
  &\text{if $\<\u,\v,\w\>$ is a permutation of a triple in $\V(\gc\B)$,}
    \\\notag&\text{then }\g(\u,\v)\rp\g(\v,\w)\geq\g(\u,\w).
\end{align*}
Define a binary relation $\Xi(\u)$ for each $\u\in\U(\gc\B)$ as
follows.  If $|\u|=1$ then $\Xi(\u)=\emptyset$.  If $|\u|>1$ then
$\Xi(\u)=\{\<\u,(\p\l)^{\gc\B}\>,\<\u,(\p\m)^{\gc\B}\>\}$ where $\p$
is the unique reduced trail in $\u$, hence $\u=\p^\gc\B$, $\k$ is the
pointer of $\p$, \ie, $\p=\p\k$, and $\{0,1,2\}=\{\k,\l,\m\}$.  In
this case, $ \Xi(\u)=\{\<\u,(\p\lambda)^{\gc\B}\>:
\text{$\p=\p\k\in\u$ is reduced},\k\neq\l<3\}.$ From the definitions
it follows that $\Xi(\u)\subseteq\V_2$ and
$(\Xi(\u))^{-1}\subseteq\V_2$, since the ordered pairs in $\Xi(\u)$
occur in the triple of the reduced trail in $\u$, which is a
permutation of $\<\u,(\p\l)^{\gc\B},(\p\m)^{\gc\B}\>$.  Points
$(\p\lambda)^{\gc\B}$ and $(\p\mu)^{\gc\B}$ (which may coincide) have
shorter reduced trails (and hence are distinct from $\u$) because
condition (ii) in the definition of $\approx$ applies to $\p\lambda$
and $\p\mu$ (due to $\kappa\neq\lambda,\mu$).  In computing their
reduced trails, only conditions (ii) and (iii) in the definition of
$\approx$ are used.  Condition (i) cannot apply because $\p$ is
reduced.  Thus, if $\<\v,\w\>\in\Xi(\u)$, either $\v$ has a shorter
reduced trail than $\w$, \ie, $|\v|<|\w|$, or $\w$ has a shorter
reduced trail than $\v$, \ie, $|\w|<|\v|$. We now turn to proving that
$\V_2$ can be partitioned into three pieces, the set $Id_{\U(\gc\B)}$
of {\bf identity pairs}, the set $\Xi(\U)$ of {\bf $\Xi$-pairs}, and
the set $\B_2$ of {\bf base-pairs}:
\begin{align*}
  \V_2&=Id_{\U(\gc\B)}\cup\Xi(\U)\cup\B_2,\\
  \Xi(\U)&=\bigcup_{\u\in\U(\gc\B)}\(\Xi(\u)\cup(\Xi(\u))^{-1}\),\\
  \B_2&=\{\<(\<\t,\kappa\>)^{\gc\B},(\<\t,\lambda\>)^{\gc\B}\>:\t\in\B,
  \,\t_\m\leq\di,\, \{\k,\l,\m\}=\{0,1,2\}\}.
\end{align*}
A pair $\<\u,\v\>$ with $\u,\v\in\U(\gc\B)$ is an identity pair if
$\u=\v$, a $\Xi$-pair if $\<\u,\v\>\in\Xi(\u)$ or
$\<\v,\u\>\in\Xi(\v)$, and a base-pair if it occurs in the triple of a
trail $\p$ of length 1.  In connection with the following lemma, it is
worth observing that $\V(\gc\B)$ is not closed under the permutation
of triples.
\begin{lemma}\label{lem8}
  If $\<\u_0,\u_1,\u_2\>\in\V(\gc\B)$ then
  $\<\u_{\pi(0)},\u_{\pi(1)},\u_{\pi(2)}\>\in\V(\gc\B)$ for every
  non-permutation $\pi\colon\{0,1,2\}\to\{0,1,2\}$.
\end{lemma}
\proof Assume $\<\u_0,\u_1,\u_2\>\in\V(\gc\B)$.  Then there is a
reduced trail
\begin{equation*}
  \p=\<\t_0,\kappa_0,\cdots,\t_\n,\kappa_\n\>\in\Tr(\gc\B),
\end{equation*}
where $\t_\n=\a=\<\a_0,\a_1,\a_2\>\in\B$ and $\a_0,\a_1,\a_2\in\atm$,
such that
$\<\u_0,\u_1,\u_2\>=\<(\p0)^{\gc\B},(\p1)^{\gc\B},(\p2)^{\gc\B}\>.$
For distinct $\k,\l<3$, define the non-permutation (replacement)
$[\k/\l]\colon3\to3$ by $[\k/\l]=\{\<k,\l\>,\<\l,\l\>,\<\m,\m\>\}$
where $\{\k,\l,\m\}=\{0,1,2\}$.  Let
\begin{align*}
  \p_\pi=\begin{cases}
    \p\odot\<\a,0,\<\a_0,\a_0,\dom{\a_0}\>,0\>&\text{ if }\pi=[1/0],\\
    \p\odot\<\a,0,\<\a_0,\rng{\a_0},\breve{\a_0}\>,0\>&\text{ if }\pi=[2/0],\\
    \p\odot\<\a,1,\<\a_1,\a_1,\dom{\a_1}\>,0\>&\text{ if }\pi=[0/1],\\
    \p\odot\<\a,1,\<\rng{\a_1},\a_1,\a_1\>,0\>&\text{ if }\pi=[2/1],\\
    \p\odot\<\a,2,\<\breve{\a_2},\dom{\a_2},\a_2\>,0\>&\text{ if }\pi=[0/2],\\
    \p\odot\<\a,2,\<\rng{\a_2},\a_2,\a_2\>,0\>&\text{ if }\pi=[2/2].
\end{cases}
\end{align*}
Then, as simple computational checks will show,
\begin{equation*}
\<(\p_\pi0)^{\gc\B},(\p_\pi1)^{\gc\B},(\p_\pi2)^{\gc\B}\>=
\<(\p\pi(0))^{\gc\B},(\p\pi(1))^{\gc\B},(\p\pi(2))^{\gc\B}\>\in\V(\gc\B).
\end{equation*}
Thus the lemma holds for the six listed replacements. The remaining 15
non-permutations can be obtained from these by composition in multiple
ways, so the lemma holds for all of them.
\endproof
\begin{lemma}\label{lem9}
  Every triple in $\V(\gc\B)$ is a permutation of one of the following
  six triples, for some distinct $\u,\v,\w\in\U(\gc\B)$.
  \begin{enumerate}
  \item[\rm(i)] $\<\u,\v,\w\>$ where $\{\<\u,\v\>,\<\u,\w\>\}=\Xi(\u)$ and
    $\<\v,\w\>$ is a $\Xi$-pair,
  \item[\rm(ii)] $\<\u,\v,\w\>$ where $\{\<\u,\v\>,\<\u,\w\>\}=\Xi(\u)$, and
    $\<\v,\w\>$ is a base-pair,
  \item[\rm(iii)] $\<\u,\v,\w\>$ where $\<\u,\v\>$, $\<\u,\w\>$, and $\<\v,\w\>$
    are base-pairs,
  \item[\rm(iv)] $\<\u,\u,\v\>$ where $\<\u,\v\>$ is a $\Xi$-pair,
  \item[\rm(v)] $\<\u,\u,\v\>$ where $\<\u,\v\>$ is a base-pair,
  \item[\rm(vi)] $\<\u,\u,\u\>$.
  \end{enumerate}
\end{lemma}
\proof Triples of types (i)--(iii) contain three distinct points in
$\U(\gc\B)$, triples of types (iv) and (v) have two distinct points,
while triples of type (vi) contain just a single point (repeated
twice). Types (i), (ii), and (iv) are the triples of trails that have
length $>1$, while types (iii) and (v) are the triples of trails of
length $1$. Triples of type (vi) are the triples of trails of all
lengths.  We proceed by induction on the length of trails.

Start with a trail $\p=\<\t,\nu\>$ of length 1, where $\t\in\B$, and
$\nu<3$.  If $\t$ is a diversity cycle, then $\t$ is not in $\E_{01}$,
$\E_{02}$, or $\E_{12}$ since $\t_2\not\leq\id$, $\t_1\not\leq\id$,
and $\t_0\not\leq\id$, respectively.  The trails $\p0$, $\p1$, and
$\p2$ are reduced, for in the definition of reduced trail, the
hypothesis of condition (i) is not met because
$\t\notin\E_{01}\cup\E_{02}\cup\E_{12}$, and the hypotheses of
condition (ii) and (iii) are not met because $|\p|=1$.  The three
reduced trails $\p0$, $\p1$, and $\p2$ are distinct, because they have
distinct pointers, so their equivalence classes (points in
$\U(\gc\B)$) are distinct, and the triple
$\<(\p0)^{\gc\B},(\p1)^{\gc\B},(\p2)^{\gc\B}\>$ of $\p$ has
type (iii), \ie, $(\p0)^{\gc\B}$, $(\p1)^{\gc\B}$, and $(\p2)^{\gc\B}$
are distinct.

Assume $\t$ is an identity cycle, \ie, $\t$ has the form
$\<\a,\a,\dom\a\>$, $\<\rng\a,\a,\a\>$, $\<\a,\rng\a,\con\a\>$, or
$\<\dom\a,\dom\a,\dom\a\>$, for some atom $\a\in\atm.$ If
$\t=\<\a,\a,\dom\a\>$, then $\t\in\E_{01},$ $\p0=\<\t,0\>$ and
$\p2=\<\t,1\>$ are reduced so $(\p0)^{\gc\B}$ and $(\p2)^{\gc\B}$ are
distinct, $\p0=\<\t,0\>\approx\<\t,1\>=\p1$ since $\t_2\leq\id$, so
$(\p0)^{\gc\B}=(\p1)^{\gc\B}$, and the triple of $\p$ has the form
$\<\u,\u,\v\>$, which is type (v).  Similarly, if
$\t=\<\rng\a,\a,\a\>$, then $\t\in\E_{12},$ $\p0$ and $\p1$ are
reduced, $\p1\approx\p2$, and the triple of $\p$ has the form
$\<\v,\u,\u\>$, which is type (v).  If $\t=\<\a,\rng\a,\con\a\>$ then
$\t\in\E_{02},$ $\p0$ and $\p1$ are reduced, $\p0\approx\p2$, and the
triple of $\p$ has the form $\<\u,\v,\u\>$, again of type (v).
Finally, if $\t=\<\dom\a,\dom\a,\dom\a\>$ then
$\t\in\E_{01}\cup\E_{12}$, $\p0$ is reduced,
$\p0\approx\p1\approx\p2$, and the triple of $\p$ has type (vi).

As inductive hypothesis, we assume that the triple of every trail of
length $\n$ or less is exactly one of the six types.  Let $\p$ be a
trail of length $\n+1$, \ie, there is a trail $\q\in\Tr(\gc\B)$ such
that
\begin{align*}
  \q&=\<\t_0,\kappa_0,\cdots,\t_{\n-1},\kappa_{\n-1}\>,\\
  \p&=\<\t_0,\kappa_0,\cdots,\t_{\n-1},\kappa_{\n-1},\t_\n,\kappa_\n\>.
\end{align*}
We may assume $\p$ is not subject to any reductions of type (i)
because if it were, it would be $\approx$ to a shorter trail whose
triple is the same as that of $\p$, hence subject to the inductive
hypothesis, which yields the desired conclusion. Therefore we assume
$\p$ is (i)-reduced.  Let $\k=\k_{\n-1}$ and $\{\k,\l,\m\}=\{0,1,2\}$.
Note that $\p\l$ and $\p\m$ are not reduced, and indeed
$(\q\l)^{\gc\B}=(\p\l)^{\gc\B}$ and $(\q\m)^{\gc\B}=(\p\m)^{\gc\B}$
since, by a (ii)-reduction,
\begin{align*}
  \p\l&=\<\t_0,\kappa_0,\cdots,\t_{\n-1},\k_{\n-1},\t_\n,\l\>\approx
  \<\t_0,\kappa_0,\cdots,\t_{\n-1},\l\>=\q\l,\\
  \p\m&=\<\t_0,\kappa_0,\cdots,\t_{\n-1},\k_{\n-1},\t_\n,\m\>\approx
  \<\t_0,\kappa_0,\cdots,\t_{\n-1},\m\>=\q\m.
\end{align*}
The inductive hypothesis applies to the trail $\q$ of length $\n$.

Suppose $\t_\n\notin\E_{\k\l}\cup\E_{\k\m}$.  This implies that $\p\k$
is reduced, because $\p$ is (i)-reduced and no reduction of types (ii)
or (iii) can apply to $\p\k$. Hence, by definition, $
\Xi((\p\k)^{\gc\B})=\{\<(\p\k)^{\gc\B},(\p\l)^{\gc\B}\>,
\<(\p\k)^{\gc\B},(\p\m)^{\gc\B}\>\}.$ If
$(\q\l)^{\gc\B}=(\q\m)^{\gc\B}$ then the triple of $\p$ is the second
kind of type (iv). Assume $(\q\l)^{\gc\B}\neq(\q\m)^{\gc\B}$.
Depending on the type of the triple of $\q$, the diversity pair
$\<(\q\l)^{\gc\B},(\q\m)^{\gc\B}\>$ is either a $\Xi$-pair, which
means the triple of $\p$ is type (i), or a base-pair and the triple of
$\p$ is type (ii).  The triple of $\p$ cannot be type (vi).

Suppose $\t_\n\in\E_{\k\l}$.  From this we get $\p\k\approx\p\l$ by a
reduction of type (iii), so $(\p\k)^{\gc\B}=(\p\l)^{\gc\B}$. Recall
that $(\q\l)^{\gc\B}=(\p\l)^{\gc\B}$ and
$(\q\m)^{\gc\B}=(\p\m)^{\gc\B}$, so if $(\q\l)^{\gc\B}=(\q\m)^{\gc\B}$
then the triple of $\p$ is type (vi). Assume
$(\q\l)^{\gc\B}\neq(\q\m)^{\gc\B}$.  From the inductive hypothesis
applied to $\q$ we know the diversity pair
$\<(\q\l)^{\gc\B},(\q\m)^{\gc\B}\>$ is either a $\Xi$-pair, in which
case the triple of $\p$ is type (iv), or it is a base-pair and the
triple of $\p$ is type (v).  Similarly, if $\t_\n\in\E_{\k\m}$ then
the triple of $\p$ is type (iv) or (v).
\endproof
The following result is the analogue of \cite[Theorem~C]{MR987611} for
weakly associative relation algebras.
\begin{theorem}\label{th7}
  Suppose $\gc\B$ is the suitable structure of a complete atomic
  $\gc\A\in\WA$.  Then $\gc\A\cong\Rl_{\S(\gc\B)}\gc\A'$ for some set
  relation algebra $\gc\A'\subseteq\Re\(\U(\gc\B)\)$, where
\begin{align*}
  \S(\gc\B)&=\{\<\u,\v\>:\u,\v\in\U(\gc\B),\,
  \exists\w\big(\<\u,\v,\w\>\in\V(\gc\B)\big)\},\\
  \U(\gc\B)&=\lb\p^{\gc\B}:\p\in\Tr(\gc\B)\rb,\\
  \V(\gc\B)&=\bigcup_{\t\in\B}\rep\t,\\
  \rep\t&=\lb\left<(\p\k)^{\gc\B}:\k<3\right>:\p\in\Tr(\gc\B),
  \ \p\text{ ends at }\t\rb,\\
  \p^{\gc\B}&=\lb\p':\p\,\approx\,\p'\rb.
\end{align*}
\end{theorem}
\proof For brevity, use ``$\Cyl_\k$'' and ``$\Diag_{\k\l}$'' in place
of ``$\Cyl_\k^{[\un]}$'' and ``$\Diag_{\k\l}^{[\un]}$''.  Let
$\E(\gc\B)=\Cyl_2\V(\gc\B)$.  Then, as we will show,
$\Ra\Rl_{\V(\gc\B)}\(\gc\C\)\cong\Rl_{\E(\gc\B)}\(\Ra\gc\C\)$ via the
isomorphism $\X\mapsto\Cyl_2\X$, where $\X$ is an element of
$\Ra\Rl_{\V(\gc\B)}\(\gc\C\)$, \ie, $\V(\gc\B)\cap\Cyl_2\X=\X\in\C.$
Indeed, if $\Y=\Cyl_2\X$, then $\X$ can be recovered from $\Y$ since
$\X=\V(\gc\B)\cap\Y$, and $\Y$ is an element of
$\Rl_{\E(\gc\B)}\(\Ra\gc\C\)$, \ie,
$\E(\gc\B)\supseteq\Cyl_2\Y=\Y\in\C,$ because $\Y\in\C$ since $\gc\C$
is a subalgebra (closed under cylindrifications), $\Cyl_2\Y=\Y$ since
cylindrifying twice is the same as doing it once, and
$\Y\subseteq\E(\gc\B)$ since $\X\subseteq\V(\gc\B)$, by the
monotonicity of cylindrification. The map $\X\mapsto\Cyl_2\X$ is thus
a bijection. We need to show it preserves relative multiplication and
converse. Toward this end, suppose $\X,\Y$ are elements of
$\Ra\Rl_{\V(\gc\B)}\(\gc\C\)$, that is,
$\V(\gc\B)\cap\Cyl_2\X=\X\in\C$ and $\V(\gc\B)\cap\Cyl_2\Y=\Y\in\C.$
For relative product we will prove
\begin{align*}\tag{$*{**}$}
  &\E(\gc\B)\cap\Cyl_2\big(\Cyl_1(\Diag_{12}\cap\Cyl_2\X)\cap\Cyl_0
  (\Diag_{02}\cap\Cyl_2\Y)\big)\\
  &=\Cyl_2\Big(\V(\gc\B)\cap\Cyl_1(\Diag_{12}\cap\X)
  \cap\Cyl_0(\Diag_{02}\cap\Y)\Big).
\end{align*}
The left side of ($*{**}$) is obtained by applying the map $\Cyl_2$ to
$\X$ and $\Y$, and then computing the relative product of their images
in the target algebra $\Rl_{\E(\gc\B)}\(\Ra\gc\C\)$, while the right
side is the result of applying the map $\Cyl_2$ to the relative
product of $\X$ and $\Y$, as computed in the source algebra
$\Ra\Rl_{\V(\gc\B)}\(\gc\C\)$, and simplifying using
$\cyl_2\x\cdot\cyl_2\y=\cyl_2(\x\cdot\cyl_2\y)$.  The right side of
($*{**}$) is included in the left side, by just the monotonicity and
idempotence of cylindrification, so we need only show the left side is
included in the right side.  Suppose $\<\u_0,\u_1,\u_2\>$ is in the
left side of (${*}{**}$), \ie,
\begin{align*}
  \<\u_0,\u_1,\u_2\>&\in\E(\gc\B)\cap\Cyl_2\(\Cyl_1(\Diag_{12}\cap\Cyl_2\X)
  \cap\Cyl_0(\Diag_{02}\cap\Cyl_2\Y)\).
\end{align*}
Then $\<\u_0,\u_1,\u_2\>\in\E(\gc\B)=\Cyl_2\(\V(\gc\B)\)$, so by the
definitions of cylindrification and $\V(\gc\B)$, there is some trail
$\p$ such that 
\begin{align*}
  \<(\p0)^{\gc\B},(\p1)^{\gc\B},(\p2)^{\gc\B}\>=\<\u_0,\u_1,\u_2\>
  &\in\Cyl_2\(\Cyl_1(\Diag_{12}
  \cap\Cyl_2\X)\cap\Cyl_0(\Diag_{02}\cap\Cyl_2\Y)\),
\end{align*}
so, by the definitions of cylindrification and diagonal elements,
there is some $\v\in\U(\gc\B)$ such that
\begin{align*}\tag{\#}
  \<\u_0,\v,\v\>&\in\Cyl_2\X,&\<\v,\u_1,\v\>&\in\Cyl_2\Y.\\
  \intertext{By the definition of cylindrification, there exist
    $\x,\y\in\U(\gc\B)$ such that} \<\u_0,\v,\x\>&\in\X,
  &\<\v,\u_1,\y\>&\in\Y.\\
  \intertext{We assumed $\X,\Y\subseteq\V(\gc\B)$, so there are trails
    $\q,\r$ such that}
  \<\u_0,\v,\x\>&=\<(\q0)^{\gc\B},(\q1)^{\gc\B},(\q2)^{\gc\B}\>,
  &\<\v,\u_1,\y\>&=\<(\r0)^{\gc\B},(\r1)^{\gc\B},(\r2)^{\gc\B}\>.
  \intertext{By Lemma \ref{lem8} we have
    $\<\u_0,\v,\v\>,\<\v,\u_1,\v\>\in\V(\gc\B)$, so by (\#) and our
    assumption that $\V(\gc\B)\cap\Cyl_2\X=\X$ and
    $\V(\gc\B)\cap\Cyl_2\Y=\Y$, we get}\tag{\#\#}
  \<\u_0,\v,\v\>&\in\Diag_{12}\cap\X,
  &\<\v,\u_1,\v\>&\in\Diag_{02}\cap\Y,
\end{align*}
hence
\begin{align*}
  \<\u_0,\u_1,\v\>&\in\Cyl_1(\Diag_{12}\cap\X)
  \cap\Cyl_0(\Diag_{02}\cap\Y).
\end{align*}
Note that if $\<\u_0,\u_1,\v\>$ happens to also be in $\V(\gc\B)$,
then $\<\u_0,\u_1,\u_2\>$ will be in $\Cyl_2\{\<\u_0,\u_1,\v\>\}$, and
will consequently be in the right hand side of (${**}*$).  We can
prove $\<\u_0,\u_1,\v\>\in\V(\gc\B)$ if $\u_0,\u_1,\v$ are not
distinct.  Since $\<\u_0,\u_1\>$ occurs in the triple of $\p$, we get
$\<\u_0,u_1,u_0\>,\<\u_0,\u_1,\u_1\>\in\V(\gc\B)$ by Lemma \ref{lem8}.
Consequently, if $\u_0=\v$ then
$\<\u_0,\u_1,\v\>=\<\u_0,u_1,u_0\>\in\V(\gc\B)$, and if $\u_1=\v$ then
$\<\u_0,\u_1,\v\>=\<\u_0,\u_1,\u_1\>\in\V(\gc\B)$.  If $\u_0=\u_1$
then $\<\u_0,\u_1,\v\>=\<\u_0,\u_0,\v\>
=\<(\q0)^{\gc\B},(\q0)^{\gc\B},(\q1)^{\gc\B}\>\in\V(\gc\B)$ by
Lemma \ref{lem8} ($\r$ can also be used here).  We may therefore
assume that $\u_0,\u_1,\v$ are distinct.  By Lemma \ref{lem9} each of
the pairs $\<\u_0,\u_1\>$, $\<\u_0,\v\>$, and $\<\v,\u_1\>$ must be
either a $\Xi$-pair or a base-pair.

Let $\w\in\{\u_0,\u_1,\v\}$ be any point with maximum length, \ie,
$|\w|\geq|\u_0|,|\u_1|,|\v|$, and let $\{\w,\x,\y\}=\{\u_0,\u_1,\v\}$.
We will show that the triple of the reduced trail in $\w$ is a
permutation of $\<\u_0,\u_1,\v\>$.  There are two cases.

Case 1: Some point in $\{\u_0,\u_1,\v\}$ has length $>1$.  By our
assumption, $|\w|>1$.  Since endpoints of base-pairs have length 1,
$\<\w,\x\>$ must be a $\Xi$-pair, hence either $\<\w,\x\>\in\Xi(\w)$
and $|w|>|x|$, or else $\<\x,\w\>\in\Xi(\x)$ and $|x|>|w|$.  But
$|x|\not>|w|$ by the choice of $\w$, so $\<\w,\x\>\in\Xi(\w)$, and,
similarly, $\<\w,\y\>\in\Xi(\w)$.  Since $\x\neq\y$, $\Xi(\w)$ has at
least two elements. By its definition, $\Xi(\w)$ has at most two
elements, so $\Xi(\w)=\{\<\w,\x\>,\<\w,\y\>\}$. Since the ordered
pairs in $\Xi(\w)$ occur in the triple of the reduced trail in $\w$,
the triple of the reduced trail in $\w$ must be a permutation of
$\<\u_0,\u_1,\v\>$.  

Case 2: $|\u_0|=|\u_1|=|\v|=1$. In this case $\w$ can be any one of
the three points. By Lemma \ref{lem9}, $\<\u_0,\u_1\>$, $\<\u_0,\v\>$,
and $\<\u_1,\v\>$ are base-pairs, \ie, the triple has type (iii) in
that lemma.  (Types (i) and (ii) require one of the points to have
length $>1$.  Types (iv), (v), and (vi) have non-distinct triples.)
Recall that the points arise from trails $\p,\q,\r$ by
$\u_0=(\p0)^{\gc\B}=(\q0)^{\gc\B}$,
$\u_1=(\p1)^{\gc\B}=(\r1)^{\gc\B}$, and
$\v=(\q1)^{\gc\B}=(\r0)^{\gc\B}$, which imply $\p0\approx\q0$,
$\p1\approx\r1$, and $\q1\approx\r0$.  Therefore $\p,\q,\r$ all begin
at the same cycle, say $\t\in\B$, because $\approx$-reductions only
relate trails that begin at the same cycle.  In general, the reduced
trail of a point $\u$ is always a ``subtrail'' of any trail that has
$\u$ in its triple. Type (i) reductions extract parts in the middle of
a trail. Type (ii) reductions just change the pointer.  Type (iii)
reductions shorten a trail by deleting the end and the ordinal
preceding it.  The reduced trail of a point $\u$ of length 1 has
length 1 (is of the form $\<\t,\k\>$) and its beginning $\t$ is the
beginning of every trail whose triple includes $\u$.  In the current
case, the trails $\p,\q,\r$ and the indices $0,1$ produce three
distinct points, the most a trail can produce.  Since the three points
$\u_0,\u_1,\v$ have length 1, their reduced trails have length 1, \ie,
have the same cycle, say $\t\in\B$, as both beginning and end, and
differ only in their pointers, so to get three points we must use all
three pointers.  Therefore $\{\u_0,\u_1,\v\}= \{(\<\t,0\>)^{\gc\B},
(\<\t,1\>)^{\gc\B}, (\<\t,2\>)^{\gc\B}\}$ and because these points are
distinct, $\t$ must be a diversity cycle, \ie, $\t_0+\t_1+\t_2
\leq\di$. Thus, in this case as well as Case 1, $\<\u_0,\u_1,\v\>$ is
a permutation of the triple of the reduced trail in $\w$.

By \thetag{K} and the definition of $\B$, we have
$\<\g(\u_1,\v),\g(\u_0,\v),\g(\u_0,\u_1)\>\in\B.$ Let
$\t=\<\t_0,\t_1,\t_2\>\in\B$ be the end of $\p$. Then
$\t_2=\g_\p(\u_0,\u_1)=\g(\u_0,\u_1)\leq\g(\u_0,\v)\rp\g(\v,\u_1),$
so, letting $\t'=\Big\<\g(\u_1,\v),\g(\u_0,\v),\g(\u_0,\u_1)\Big\>,$
we have $\t\mathrel{\T_2}\t'.$ Define a trail $\p'$ extending $\p$ by
$ \p'=\p\odot\big\<\t,2,\t',0\big\>.$ Then
$\u_0=(\p0)^{\gc\B}=(\p'0)^{\gc\B}$ and
$\u_1=(\p1)^{\gc\B}=(\p'1)^{\gc\B}$ since $\p0\approx\p'0$ and
$\p1\approx\p'1$ by type (iii) reductions. Let $\z=(\p'2)^{\gc\B}$.
Then $\<\u_0,\u_1,\z\>\in\V(\gc\B)$ because $\<\u_0,\u_1,\z\>$ is the
triple of $\p'$, hence also $\<\u_0,\z,\z\>\in\V(\gc\B)$ by Lemma
\ref{lem8}.

Since $\Rl_{\V(\gc\B)}\(\gc\C\)$ is complete and atomic, its element
$\X$ must be a join of atoms. Since each atom has the form $\rep\s$
for some $\s\in\B$ by Theorem \ref{th6}, from (\#\#) we get
$\<\u_0,\v,\v\>\in\rep\s\subseteq\X$ for some $\s\in\B$.  By the
definition of $\rep\s$, $\<\u_0,\v,\v\>$ is the triple of some trail
$\p''$ that ends at $\s$.  Therefore $\s=\Big\<\rng{\g(\u_0,\v)},
\g(\u_0,\v),\g(\u_0,\v)\Big\>.$ From $\<\u_0,\z,\z\>\in\V(\gc\B)$ it
follows that $\<\u_0,\z,\z\>$ is the triple of a trail ending at
$\Big\<\g(\z,\z),\g(\u_0,\z),\g(\u_0,\z)\Big\>.$ The triple of $\p'$
is $\<\u_0,\u_1,\z\>$, so by the relevant definitions,
\begin{align*}
  \g(\u_0,\z)&=\g_{\p'}(\u_0,\z)
  =\g_{\p'}\big((\p'0)^{\gc\B},(\p'2)^{\gc\B}\big)=\t'_1=\g(\u_0,\v),
  \\
  \g(\z,\z)&=\g_{\p'}(\z,\z)
  =\g_{\p'}\big((\p'2)^{\gc\B},(\p'2)^{\gc\B}\big)
  =\rng{\t'_1}=\g(\u_0,\v)\rng{}.
\end{align*}
It follows that $\<\u_0,\z,\z\>$ is the triple of a trail ending at
$\s$, and therefore $\<\u_0,\z,\z\>\in\rep\s\subseteq\X$. By a similar
argument, $\<\z,\u_1,\z\>\in\Y$.  Thus we have
$\<\u_0,\z,\z\>\in\Diag_{12}\cap\X$ and
$\<\z,\u_1,\z\>\in\Diag_{02}\cap\Y$, which imply, together with
$\<\u_0,\u_1,\z\>\in\V(\gc\B)$, that
$\<\u_0,\u_1,\u_2\>\in\Cyl_2\big(\{\<\u_0,\u_1,\z\>\}\big)
\subseteq\Cyl_2\Big(\V(\gc\B)\cap
\Cyl_1(\Diag_{12}\cap\X)\cap\Cyl_0(\Diag_{02}\cap\Y)\Big),$ \ie,
$\<\u_0,\u_1,\u_2\>$ is in the right hand side of ($*{**}$), as
desired. Thus relative product is preserved.  For converse, we will
prove
\begin{align*}\tag{${**}{**}$}
  &\E(\gc\B)\cap\Cyl_2(\Diag_{20}\cap\Cyl_0(\Diag_{01}\cap\Cyl_1(\Diag_{12}
  \cap\Cyl_2\X)))\\
  &=\E(\gc\B)\cap\Cyl_2(\Diag_{20}\cap
  \V(\gc\B)\cap\Cyl_0(\Diag_{01}\cap\V(\gc\B)
  \cap\Cyl_1(\Diag_{12}\cap\V(\gc\B)\cap\X))).
\end{align*}
The left side of ($*{**}*$) is obtained by applying the isomorphism
$\Cyl_2$ to $\X$, and then computing the converse of its image in the
target algebra $\Rl_{\E(\gc\B)}\(\Ra\gc\C\)$. The right side of
($*{**}*$) is obtained by computing the converse of $\X$ in the source
algebra $\Ra\Rl_{\V(\gc\B)}\(\gc\C\)$, applying the isomorphism
$\Cyl_2$, and simplifying using $\cyl_2\x\cdot\cyl_2\y =
\cyl_2(\x\cdot\cyl_2\y)$ and the definition of $\E(\gc\B)$.  The right
side of ($*{**}*$) is included in the left side, by just the
monotonicity and idempotence of cylindrification, so we need only show
the left side is included in the right side. Suppose
$\<\u_0,\u_1,\u_2\>$ is in the left side of ($*{**}*$). Then
\begin{align*}
  \<\u_0,\u_1,\u_2\>&\in
  \Cyl_2(\Diag_{20}\cap\Cyl_0(\Diag_{01}\cap\Cyl_1(\Diag_{12}
  \cap\Cyl_2\X))).
\end{align*}
By the definitions of cylindrification and diagonal elements, we
obtain successively
\begin{align*}
  \<\u_0,\u_1,\u_0\>&\in\Diag_{02}\cap\Cyl_0(\Diag_{01}\cap\Cyl_1(\Diag_{12}
  \cap\Cyl_2\X)),\\
  \<\u_1,\u_1,\u_0\>&\in\Diag_{01}\cap\Cyl_1(\Diag_{12}\cap\Cyl_2\X),\\
  \<\u_1,\u_0,\u_0\>&\in\Diag_{12}\cap\Cyl_2\X,
\end{align*}
and, finally, there is some $\x\in\U(\gc\B)$ such that
$\<\u_1,\u_0,\x\>\in\X$.  We assumed $\X\subseteq\V(\gc\B)$, so
$\<\u_1,\u_0,\x\>\in\V(\gc\B)$, hence
$\<\u_1,\u_0,\u_0\>,\<\u_1,\u_1,\u_0\>,\<\u_0,\u_1,\u_0\>\in\V(\gc\B)$
by Lemma \ref{lem8}.  Combining these facts with our assumption that
$\V(\gc\B)\cap\Cyl_2\X=\X$, we obtain successively
\begin{align*}
  \<\u_1,\u_0,\u_0\>&\in\Diag_{12}\cap\V(\gc\B)\cap\X,\\
  \<\u_1,\u_1,\u_0\>&\in\Diag_{01}\cap\V(\gc\B)\cap\Cyl_1(\Diag_{12}
  \cap\V(\gc\B)\cap\X),\\
  \<\u_0,\u_1,\u_0\>&\in\Diag_{20}\cap\V(\gc\B)\cap\Cyl_0\big(\Diag_{01}
  \cap\V(\gc\B)\cap\Cyl_1(\Diag_{12}\cap\V(\gc\B)\cap\X)\big),\\
  \<\u_0,\u_1,\u_2\>&\in\Cyl_2\Big(\Diag_{20}\cap\V(\gc\B)\cap\Cyl_0
  \big(\Diag_{01}
  \cap\V(\gc\B)\cap\Cyl_1(\Diag_{12}\cap\V(\gc\B)\cap\X)\big)\Big).
\end{align*}
Together with $\<\u_0,\u_1,\u_2\>\in\E(\gc\B)$, we conclude
$\<\u_0,\u_1,\u_2\>$ is in the right hand side of ($*{**}*$), as
desired.  This completes the proof that
$\Ra\Rl_{\V(\gc\B)}\(\gc\C\)\cong\Rl_{\E(\gc\B)}\(\Ra\gc\C\).$ Since $\gc\C$
is a subalgebra of $\Sb(\un)$, applying the $\Ra$ operator gives
$\Ra\gc\C\subseteq\Ra\Sb(\un).$ According to \cite[5.3.16(1)]{HMT85},
$\Re\(\U(\gc\B)\)\cong\Ra\Sb\(\un\)$ via the isomorphism $\f$ defined
for $\X\subseteq\U(\gc\B)\times\U(\gc\B)$ by
\begin{equation*}
  \f(\X)=\{\<\u,\v,\w\>:\u,\v,\w\in\U(\gc\B),\,\<\u,\v\>\in\X\}
  =\X\times\U(\gc\B).
\end{equation*}
Hence there is an algebra $\gc\A'\subseteq\Re(\U(\gc\B))$ such that
$\Ra\gc\C\cong\gc\A'$ via the isomorphism $\f^{-1}$, where, if
$\Y\subseteq\un$, then
\begin{equation*}
  \f^{-1}(\Y)=\{\<\u,\v\>:\u,\v\in\U(\gc\B),\,
  \exists\w\big(\<\u,\v,\w\>\in\Y\big)\}.
\end{equation*}
Let $\S(\gc\B)=\f^{-1}(\E(\gc\B))$. Then it is easy to prove
\begin{equation*}
  \S(\gc\B)=\{\<\u,\v\>:\u,\v\in\U(\gc\B),\,
  \exists\w\big(\<\u,\v,\w\>\in\V(\gc\B)\big)\}.
\end{equation*}
Relativize $\Ra\gc\C$ and its image $\gc\A'$ under the isomorphism
$\f^{-1}$ to the ternary relation $\E(\gc\B)$ and its image
$\S(\gc\B)$ under $\f^{-1}$, obtaining
$\Rl_{\E(\gc\B)}\(\Ra\gc\C\)\cong\Rl_{\S(\gc\B)}\(\gc\A'\).$ It was
shown above that
$\Ra\Rl_{\V(\gc\B)}\(\gc\C\)\cong\Rl_{\E(\gc\B)}\(\Ra\gc\C\),$ and
from Theorem \ref{th6} we have
$\gc\A\cong\Ra\Rl_{\V(\gc\B)}\(\gc\C\),$ so we conclude that
$\gc\A\cong\Rl_{\S(\gc\B)}\(\gc\A'\).$
\endproof

\section{Elementary Laws of $\WA$}\label{sect11}
Throughout this section, we assume $\gc\A$ is an arbitrary algebra
satisfying axioms \eqref{ra1}--\eqref{ra10}.  Few specific references
are made to elementary facts from the theory of Boolean algebras, so
there are no references to axioms \eqref{ra1}--\eqref{ra3}. Axiom
\eqref{ra4} is not needed until \eqref{zx1880}.
\begin{equation}\label{1'sym}
  \con\id=\id
\end{equation}
\proof
\begin{align*}
  &\con\id\overset{\eqref{ra6}}=\con\id\rp\id
  \overset{\eqref{ra7}}=\con\id\rp\con{\con\id{}}
  \overset{\eqref{ra9}}=\big(\con\id\rp\id\big)\,\con{}
  \overset{\eqref{ra6}}=\con{\con\id}\overset{\eqref{ra7}}=\id.
\end{align*}
\endproof
\begin{equation}\label{1sym}
  \con1=1
\end{equation}
\proof
\begin{align*}
  &1=1+\con1 \overset{\eqref{ra7}}=\con{\con1}+\con1
  \overset{\eqref{ra8}}=\con{\con1+1}=\con1.
\end{align*}
\endproof
\begin{equation}\label{0sym} 
  \con0=0
\end{equation}
\proof
\begin{align*}
  &\con0=0+\con0\overset{\eqref{ra7}}=\con{\con0}+\con0
  \overset{\eqref{ra8}}=\con{\con0+0}=\con{\con0}\overset{\eqref{ra7}}=0.
\end{align*}
\endproof
\begin{equation}\label{conmon}
  \x\leq\y\implies\con\x\leq\con\y
\end{equation}
\proof\begin{align*} \x\leq\y&\iff\y=\x+\y\implies
\con\y=\con{\x+\y}\overset{\eqref{ra8}}
=\con\x+\con\y\iff\con\x\leq\con\y.
\end{align*}\endproof
\begin{equation}\label{rightmon}
  \x\leq\y\implies\x\rp\z\leq\y\rp\z
\end{equation}
\proof\begin{align*}
  \x\leq\y&\iff\y=\x+\y\implies
  \y\rp\z=(\x+\y)\rp\z\overset{\eqref{ra5}}=\x\rp\z+\y\rp\z
  \iff\x\rp\z\leq\y\rp\z.
\end{align*}\endproof
\begin{equation}\label{dual5}
  \z\rp(\x+\y)=\z\rp\x+\z\rp\y
\end{equation}
\proof
\begin{align*}
  \z\rp(\x+\y) &\overset{\eqref{ra7}}=\con{\con{\z\rp(\x+\y)}}
  \overset{\eqref{ra9}}=\con{\con{\x+\y}\rp\con\z}
  \overset{\eqref{ra8}}=\con{(\con\x+\con\y)\rp\con\z}
  \\ &\overset{\eqref{ra5}}=\con{\con\x\rp\con\z+\con\y\rp\con\z}
  \overset{\eqref{ra9}}=\con{\con{\z\rp\x}+\con{\z\rp\y}}
  \\ &\overset{\eqref{ra8}}=\con{\con{\z\rp\x+\z\rp\y}}
  \overset{\eqref{ra7}}=\z\rp\x+\z\rp\y.
\end{align*}
\endproof
\begin{equation}\label{leftmon}
  \x\leq\y\implies\z\rp\x\leq\z\rp\y
\end{equation}
\proof The proof is similar to that of \eqref{rightmon}, using
\eqref{dual5} instead of \eqref{ra5}.  \endproof
\begin{equation}\label{x;0}
 \x\rp0=0
\end{equation}
\proof
\begin{align*}
  \x\rp0 \overset{\eqref{rightmon}}\leq1\rp0
  \overset{\eqref{leftmon}}\leq1\rp\min{1\rp1}
  \overset{\eqref{1sym}}=\con1\rp\min{1\rp1}
  \overset{\eqref{ra10}}\leq\min1=0.
\end{align*}
\endproof
\begin{equation}\label{0;x}
  0\rp\x=0
\end{equation}
\proof
\begin{align*}
  0\rp\x \overset{\eqref{0sym}}=\con0\rp\x
  \overset{\eqref{ra7}}=\con0\rp\con{\con\x}
  \overset{\eqref{ra9}}=\con{\con\x\rp0} \overset{\eqref{x;0}}=\con0
  \overset{\eqref{0sym}}=0.
\end{align*}
\endproof
\begin{equation}\label{leftid}
  \id\rp\x=\x
\end{equation}
\proof
\begin{align*}
  \id\rp\x\overset{\eqref{1'sym}}=\con\id\rp\x
  \overset{\eqref{ra7}}=\con\id\rp\con{\con\x}
  \overset{\eqref{ra9}}=\con{\con\x\rp\id}
  \overset{\eqref{ra6}}=\con{\con\x} \overset{\eqref{ra7}}=\x.
\end{align*}
\endproof
\begin{equation}
  \label{zx1766}\x\rp\y\bp\z=\x\rp(\y\bp\con\x\rp\z)\bp\z
  =(\x\bp\z\rp\con\y)\rp\y\bp\z
\end{equation}
\proof We prove only the first equality in \eqref{zx1766}. The proof
of the second is quite similar.
\begin{align*}
  \x\rp\y\bp\z&=\x\rp(\y\bp\con\x\rp\z+\y\bp\min{\con\x\rp\z})\bp\z
  \overset{\eqref{dual5}}=\Big(\x\rp(\y\bp\con\x\rp\z)+
  \x\rp(\y\bp\min{\con\x\rp\z})\Big)\bp\z\\
  &\overset{\eqref{ra7}}\leq\Big(\x\rp(\y\bp\con\x\rp\z)+
  \con{\con\x}\rp(\min{\con\x\rp\z})+\min\z\Big)\bp\z
  \overset{\eqref{ra10}}=\Big(\x\rp(\y\bp\con\x\rp\z)+
  \min\z\Big)\bp\z\\
  &=\x\rp(\y\bp\con\x\rp\z)\bp\z\overset{\eqref{leftmon}}
  \leq\x\rp\y\bp\z.
\end{align*}
\endproof
\begin{align}
  \label{zx1761}&\x\rp\y\bp\z\leq\x\rp(\y\bp\con\x\rp\z)\\
  \label{zx1762}&\y\rp\x\bp\z\leq(\y\bp\z\rp\con\x)\rp\x
\end{align}
\proof Propositions \eqref{zx1761} and \eqref{zx1762} follow by
Boolean algebra from \eqref{zx1766}, and have the following
consequences, because of \eqref{0;x} and \eqref{x;0}.
\begin{align}
  \label{cycle1}&\x\rp\y\bp\z\neq0\implies\y\bp\con\x\rp\z\neq0\\
  \label{cycle2}&\y\rp\x\bp\z\neq0\implies\y\bp\z\rp\con\x\neq0
\end{align}
\endproof
\begin{equation}\label{zx1850}
  \u,\v\leq\id\implies\u\cdot\v=\u\rp\v
\end{equation}
\proof
\begin{align*}
  \u\bp\v&\overset{\eqref{ra6}}=\u\rp\id\bp\v
  \overset{\eqref{zx1766}}=\u\rp(\id\bp\con\u\rp\v)\bp\v
  \overset{\eqref{leftmon}}\leq\u\rp(\con\u\rp\v)
  \overset{\eqref{conmon}\eqref{rightmon}\eqref{leftmon}}
  \leq\u\rp(\con\id\rp\v)\\
  &\overset{\eqref{1'sym}}=\u\rp(\id\rp\v)
  \overset{\eqref{leftid}}=\u\rp\v
  \overset{\eqref{rightmon}\eqref{leftmon}}\leq\id\rp\v\bp\u\rp\id
  \overset{\eqref{ra6}\eqref{leftid}}=\u\bp\v.
\end{align*}
\endproof
\begin{equation}\label{zx1849a}
  \u\leq\id\implies\u\leq\con\u
\end{equation}
\proof
\begin{align*}
  \u\overset{\eqref{ra6}}=\u\rp\id\bp\u
  \overset{\eqref{zx1761}}\leq\u\rp(\id\bp\con\u\rp\u)
  \overset{\eqref{leftmon}}\leq\u\rp(\con\u\rp\u)
  \overset{\eqref{rightmon}}\leq\id\rp(\con\u\rp\u)
  \overset{\eqref{leftid}}=\con\u\rp\u
  \overset{\eqref{leftmon}}\leq\con\u\rp\id
  \overset{\eqref{ra6}}=\con\u.
\end{align*}
\endproof
\begin{equation}\label{zx1849}
  \u\leq\id\implies\con\u=\u
\end{equation}
\proof From $\u\leq\id$ get $\con\u\leq\con\id=\id$ by \eqref{conmon}
and \eqref{1'sym}. Then $\u\leq\con\u$ by \eqref{zx1849a} and
$\con\u\leq\con{\con\u}=\u$ by \eqref{zx1849a} and \eqref{ra7}, hence
we have $\con\u=\u$.\endproof
\begin{equation}\label{zx1861}
  \u\leq\id\implies\u\rp\x=\x\bp\u\rp1
\end{equation}
\proof
\begin{align*}
  \u\rp\x&\overset{\eqref{rightmon}}\leq\id\rp\x\bp\u\rp\x
  \overset{\eqref{leftid}}=\x\bp\u\rp\x
  \overset{\eqref{leftmon}}\leq\x\bp\u\rp1
  \\ &\overset{\eqref{zx1761}}\leq\u\rp(1\bp\con\u\rp\x)
  \overset{\eqref{zx1849}}=\u\rp(\u\rp\x)
  \overset{\eqref{rightmon}}\leq\id\rp(\u\rp\x)
  \overset{\eqref{leftid}}=\u\rp\x.
\end{align*}
\endproof
\begin{equation}\label{zx1699}
  \emptyset\neq\X\subseteq\A\implies\con{\sum\X}=\sum\{\con\x:\x\in\X\}
\end{equation}
\proof Assuming $\emptyset\neq\X\subseteq\A$ and $\sum\X$ exists, we
must show $\con{\sum\X}$ is the least upper bound of
$\{\con\x:\x\in\X\}$.  By \eqref{conmon}, $\con\x\leq\con{\sum\X}$ for
every $\x\in\X$ since $\x\leq\sum\X$, so $\con{\sum\X}$ is an upper
bound of $\{\con\x:\x\in\X\}$.  Assume $\con\x\leq\y$ for every
$\x\in\X$ ($\y$ is an upper bound).  Then $\x\leq\con\y$ for every
$\x\in\X$ by \eqref{ra7} and \eqref{conmon}, so $\con\y$ is an upper
bound of $\X$, hence $\sum\X\leq\con\y$.  But then
$\con{\sum\X}\leq\y$ by \eqref{ra7} and \eqref{conmon}. Thus
$\con{\sum\X}$ is the \emph{least} upper bound of
$\{\con\x:\x\in\X\}$.  \endproof
\begin{equation}\label{ca}
  \X\subseteq\A\implies\sigop{\textstyle\sum\X}
  =\textstyle\sum\{\sigop\x:\x\in\X\}
\end{equation}
\proof If $\X=\emptyset$ then $\sum\X=0$ and
$\{\sigop\x:\x\in\X\}=\emptyset$, so the result follows from
\eqref{x;0}. Assume $\emptyset\neq\X\subseteq\A$ and $\sum\X$ exists.
For any $\y\in\A$, the following statements are equivalent (mostly for
Boolean algebraic reasons).
\begin{align*} 
  &	\sigop{\textstyle\sum\X}\leq\y\\
  &	\sigop{\textstyle\sum\X}\cdot\min\y=0\\
  & \textstyle\sum\X\cdot\tauop{\min\y}=0
  &&\text{\eqref{cycle1}, \eqref{ra7}}\\
  &	\textstyle\sum\X\leq\min{\tauop{\min\y}}\\
  &	\text{$\x\leq\min{\tauop{\min\y}}$ for all $\x\in\X$}\\
  &	\text{$\x\cdot\tauop{\min\y}=0$ for all $\x\in\X$}\\
  & \text{$\sigop\x\cdot\min\y=0$ for all $\x\in\X$}
  &&\text{\eqref{cycle1}, \eqref{ra7}}\\
  &	\text{$\sigop\x\leq\y$ for all $\x\in\X$}\\
  & \text{$\y$ is an upper bound of $\{\sigop\x:\x\in\X\}$}
\end{align*}
Since the first statement is true when $\sigop{\sum\X}=\y$,
$\sigop{\sum\X}$ is itself an upper bound of $\{\sigop\x:\x\in\X\}$.
Since the last statement implies the first, $\sigop{\sum\X}$ is
included in all the upper bounds of $\{\sigop\x:\x\in\X\}$.  Thus
$\sigop{\sum\X}$ is the least upper bound of $\{\sigop\x:\x\in\X\}$,
\ie, $\sigop{\sum\X}=\sum\{\sigop\x:\x\in\X\}$. \endproof
\begin{equation}\label{1(36)}
  \x\in\atm\implies\con\x\in\atm
\end{equation}
\proof To show $\con\x$ is an atom it suffices to show that
$\con\x\neq0$ and if $0\neq\con\x\bp\y$ then $\con\x\leq\y$. Note that
if $\con\x=0$ then $\x=\con{\con\x}=\con0=0$ by \eqref{0sym} and
\eqref{ra7}. Thus $\con\x$ is not zero because $\x$ is not zero. For
the other part, we have
\begin{align*}
  0&\neq\con\x\bp\y&&\text{Hyp.}\\
  0&\neq\con\x\rp\id\bp\y&&\text{\eqref{ra6}}\\
  0&\neq\id\bp\x\rp\y&&\text{\eqref{cycle1}, \eqref{ra7}}\\
  0&\neq\id\rp\con\y\bp\x&&\text{\eqref{cycle2}}\\
  0&\neq\con\y\bp\x&&\text{\eqref{leftid}}\\
  x&\leq\con\y&&\text{$\x$ is an atom}\\
  \con\x&\leq\y&&\text{\eqref{conmon}, \eqref{ra7}}
\end{align*}
\endproof
\begin{align}\label{cyclelaw}
  \x,\y,\z\in\atm\implies \Big(&\x\rp\y\geq\z\iff\con\x\rp\z\geq\y\iff
  \y\rp\con\z\geq\con\x\iff\\\notag
  &\con\y\rp\con\x\geq\con\z\iff\con\z\rp\x\geq\con\y\iff
  \z\rp\con\y\geq\x\Big)
\end{align}
\proof Suppose $\x,\y,\z$ are atoms such that $\x\rp\y\geq\z.$ The
atom $\z$ is not zero, so $\x\rp\y\bp\z\neq0,$ from which we get
$\y\bp\con\x\rp\z\neq0$ by \eqref{cycle1}, then
$\y\rp\con\z\bp\con\x\neq0$ by \eqref{cycle2}, then
$\con\z\bp\con\y\rp\con\x\neq0$ by \eqref{cycle1}, then
$\con\z\rp\x\bp\con\y\neq0$ by \eqref{cycle2} and \eqref{ra7}, and
finally $\x\bp\z\rp\con\y\neq0$ by \eqref{cycle1} and \eqref{ra7}.
Since $\con\x,\con\y,\con\z$ are also atoms by \eqref{1(36)}, these
last five inequalities are respectively equivalent to the inequalities
$\con\x\rp\z\geq\y$, $\y\rp\con\z\geq\con\x$,
$\con\y\rp\con\x\geq\con\z$, $\con\z\rp\x\geq\con\y$, and
$\z\rp\con\y\geq\x$. \endproof
Define domain and range operators as follows.
\begin{equation}\label{domrng}
  \dom{x}=x\rp\con{x}\cdot\id\qquad\rng{x}=\con{x}\rp x\cdot\id
\end{equation}
\begin{equation}\label{condomrng}
  \dom{\con\x}=\rng\x\qquad\rng{\con\x}=\dom\x
\end{equation}
\proof
\begin{align*}
  \dom{\con\x}\overset{\eqref{domrng}}=
  \con\x\rp\con{\con\x}\bp\id\overset{\eqref{ra7}}
  =\con\x\rp\x\bp\id\overset{\eqref{domrng}}=\rng\x,\qquad
  \rng{\con\x}\overset{\eqref{domrng}}=
  \con{\con\x}\rp\con\x\bp\id\overset{\eqref{ra7}}
  =\x\rp\con\x\bp\id\overset{\eqref{domrng}}=\rng\x.
\end{align*}
\endproof
\begin{equation}\label{zx1867}
  \x=\dom\x\rp\x=\x\rp\rng\x
\end{equation}
\proof
\begin{align*}
  \x\overset{\eqref{leftid}}=\id\rp\x\bp\x
  \overset{\eqref{zx1762}}\leq(\id\bp\x\rp\con\x)\rp\x
  \overset{\eqref{domrng}}=\dom\x\rp\x
  \overset{\eqref{rightmon}}\leq\id\rp\x\overset{\eqref{leftid}}=\x,
  \\ \x\overset{\eqref{ra6}}=\x\rp\id\bp\x
  \overset{\eqref{zx1761}}\leq\x\rp(\id\bp\con\x\rp\x)
  \overset{\eqref{domrng}}=\x\rp\rng\x
  \overset{\eqref{leftmon}}\leq\x\rp\id \overset{\eqref{ra6}}=\x.
\end{align*}
\endproof
\begin{align}\label{domid}
  \u\leq\id\implies\dom\u=\rng\u=\u
\end{align}
\proof
\begin{align*}
  \u&\overset{\eqref{zx1867}}=\dom\u\rp\u
  \overset{\eqref{leftmon}}\leq\dom\u\rp\id
  \overset{\eqref{ra6}}=\dom\u
  \overset{\eqref{domrng}}=\id\bp\u\rp\con\u
  \overset{\eqref{zx1849}}=\id\bp\u\rp\u
  \overset{\eqref{leftmon}}\leq\u\rp\id \overset{\eqref{ra6}}=\u,
  \\ \u&\overset{\eqref{zx1867}}=\u\rp\rng\u
  \overset{\eqref{rightmon}}\leq\id\rp\rng\u
  \overset{\eqref{leftid}}=\rng\u
  \overset{\eqref{domrng}}=\id\bp\con\u\rp\u
  \overset{\eqref{leftmon}}\leq\con\u\rp\id
  \overset{\eqref{ra6}}=\con\u \overset{\eqref{zx1849}}=\u.
\end{align*}
\endproof
\begin{equation}\label{zx1880}
  \u\leq\id\implies(\u\rp\x)\rp\y=\u\rp(\x\rp\y)
\end{equation}
\proof
\begin{align*}
  (\u\rp\x)\rp\y &\overset{\eqref{rightmon}\eqref{leftmon}}\leq
  (\u\rp1)\rp1\bp(\id\rp\x)\rp\y
  \overset{\eqref{ra4}\eqref{leftid}}=\u\rp1\bp\x\rp\y
  \overset{\eqref{zx1861}}=\u\rp(\x\rp\y)
  \\ &\overset{\eqref{rightmon}\eqref{leftmon}}=\u\rp1\bp\id\rp(\x\rp\y)
  \overset{\eqref{leftid}}=\u\rp1\bp\x\rp\y
  \overset{\eqref{zx1762}}\leq(\x\bp(\u\rp1)\rp\con\y)\rp\y
  \\ &\overset{\eqref{rightmon}\eqref{leftmon}}
  \leq(\x\bp(\u\rp1)\rp1)\rp\y
  \overset{\eqref{ra4}}=(\x\bp\u\rp1)\rp\y
  \overset{\eqref{zx1861}}=(\u\rp\x)\rp\y.
\end{align*}
\endproof
\begin{equation}\label{zx1882}
  \u\leq\id\implies(\x\rp\u)\rp\y\leq\x\rp(\u\rp\y)
\end{equation}
\proof
\begin{align*}
  (\x\rp\u)\rp\y&=(\x\rp\u)\rp\y\bp1
  \overset{\eqref{zx1761}}\leq(\x\rp\u)\rp(y\bp\con{\x\rp\u}\rp1)
  \overset{\eqref{ra9}}\leq(\x\rp\u)\rp(y\bp(\con\u\rp\con\x)\rp1)
  \\&\overset{\eqref{zx1849}}=(\x\rp\u)\rp(y\bp(\u\rp\con\x)\rp1)
  \overset{\eqref{rightmon}\eqref{leftmon}}
  \leq(\x\rp\id)\rp(\y\bp(\u\rp1)\rp1)
  \\&\overset{\eqref{ra6}\eqref{ra4}}=\x\rp(\y\bp\u\rp1)
  \overset{\eqref{zx1861}}=\x\rp(\u\rp\y).
\end{align*}
\endproof
\begin{equation}\label{zx1883}
  \u,\v\leq\id\implies(\x\rp\u)\rp(\v\rp\y)\leq\x\rp((\u\bp\v)\rp\y)
\end{equation}
\proof
\begin{align*}
  (\x\rp\u)\rp(\v\rp\y)
  \overset{\eqref{zx1882}}\leq\x\rp(\u\rp(\v\rp\y))
  \overset{\eqref{zx1880}}=\x\rp((\u\rp\v)\rp\y)
  \overset{\eqref{zx1850}}=\x\rp((\u\bp\v)\rp\y).
\end{align*}
\endproof
\begin{equation}\label{zx1890} 
  \x\in\atm\implies\dom\x,\rng\x\in\atm
\end{equation}
\proof If $\dom\x=0$, then $\x=\dom\x\rp\x=0\rp\x=0$ by \eqref{zx1867}
and \eqref{0;x}, but $\x\neq0$ since $\x$ is an atom. Therefore
$\dom\x\neq0$. Assume $\y\in\A$. We derive $\dom\x\leq\y$ from
$\dom\x\bp\y\neq0$. First we show $\x\leq(\id\bp\y)\rp1$.
\begin{align*}
  &\dom\x\bp\y\neq0&&\text{Hyp.}\\
  &\x\rp\con\x\bp\id\bp\y\neq0&&\text{\eqref{domrng}}\\
  &\x\bp(\id\bp\y)\rp\x\neq0&&\text{\eqref{cycle2}\eqref{ra7}}\\
  &\x\leq(\id\bp\y)\rp\x\leq(\id\bp\y)\rp1 &&\text{$\x$ is an atom,
    \eqref{leftmon}}
\end{align*}
Next we use $\x\leq(\id\bp\y)\rp1$ (in the second step below) to show
$\dom\x\leq\y$. The first use of weak associativity \eqref{ra4} occurs
in the fourth step.
\begin{align*}
  \dom\x &\overset{\eqref{domrng}}=\x\rp\con\x\bp\id
  \overset{\eqref{rightmon}}\leq((\id\bp\y)\rp1)\rp\con\x\bp\id
  \overset{\eqref{leftmon}}\leq((\id\bp\y)\rp1)\rp1\bp\id\\
  &\overset{\eqref{ra4}}=(\id\bp\y)\rp1\bp\id
  \overset{\eqref{zx1761}}\leq(\id\bp\y)\rp(1\bp\con{\id\bp\y}\rp\id)
  \overset{\eqref{rightmon}}\leq\y\rp(\con{\id\bp\y}\rp\id)\\
  &\overset{\eqref{ra6}}=\y\rp\con{\id\bp\y}
  \overset{\eqref{zx1849}}=\y\rp(\id\bp\y)
  \overset{\eqref{leftmon}}\leq\y\rp\id \overset{\eqref{ra6}}=\y.
\end{align*}
\endproof
\begin{equation}\label{zx1891}
  \x,\y\in\atm\text{ and }\x\rp\y\neq0\implies\rng\x=\dom\y
\end{equation}
\proof Since $\rng\x$ and $\dom\y$ are atoms by \eqref{zx1890}, they
can fail to be equal only by being disjoint. However, if
$\rng\x\bp\dom\y=0$ then we get a contradiction as follows.
\begin{align*}
  0\overset{\text{Hyp.}}\neq\x\rp\y
  \overset{\eqref{zx1867}}=(\x\rp\rng\x)\rp(\dom\y\rp\y)
  \overset{\eqref{zx1883}}\leq\x\rp((\rng\x\bp\dom\y)\rp\y)
  \overset{\text{Hyp.}}=\x\rp(0\rp\y)
  \overset{\eqref{x;0}\eqref{0;x}}=0.
\end{align*}
\endproof
\begin{equation}\label{1(39)}
  \u,\x\in\atm\text{ and }\u\leq\id\text{ and }
  \u\rp\x\neq0\implies\u=\dom\x
\end{equation}
\proof
\begin{align*}
  \u\overset{\eqref{zx1867}}=\u\rp\rng\u
  \overset{\eqref{zx1891}}=\u\rp\dom\x\overset{\eqref{rightmon}}
  \leq\id\rp\dom\x\overset{\eqref{leftid}}=\dom\x,
\end{align*}
so $\u=\dom\x$ since $\u$ and $\dom\x$ are atoms.  \endproof
\bibliographystyle{plain}

\end{document}